\documentclass[reqno,12pt]{amsart}

\usepackage{graphics,amssymb,amsthm,amsmath}
\usepackage[english]{babel}

\usepackage[usenames,dvipsnames]{color}

 \newtheorem{Thm}{Theorem}[section]
\newtheorem{Def}[Thm]{Definition} \newtheorem{Rem}[Thm]{Remark}
\newtheorem{Lem}[Thm]{Lemma} \newtheorem{Cor}[Thm]{Corollary}
\newtheorem{Prop}[Thm]{Proposition}

\begin{document}
\fontsize{.5cm}{.5cm}\selectfont\sf

\title[Schubert Draft]{Double Quantum Schubert Cells and\\ Quantum Mutations}

\date{\today}
\author{Hans P. Jakobsen}
\address{
  Department of Mathematical Sciences\\ University of
Copenhagen\\Universitetsparken 5\\
   DK-2100, Copenhagen,
  Denmark} \email{jakobsen@math.ku.dk}
\begin{abstract}Let ${\mathfrak p}\subset {\mathfrak g}$ be a parabolic 
subalgebra of s simple finite dimensional Lie algebra over ${\mathbb C}$. To 
each pair $w^{\mathfrak a}\leq w^{\mathfrak c}$ of minimal left coset 
representatives in the quotient space $W_p\backslash W$ we construct explicitly 
a quantum seed ${\mathcal Q}_q({\mathfrak a},{\mathfrak 
c})$. We define Schubert creation and annihilation mutations and show that our 
seeds are related by such mutations. We also introduce more elaborate seeds to 
accommodate our mutations. The quantized Schubert Cell decomposition 
of the quantized generalized flag manifold can be viewed as the result of such 
mutations having their origins in the pair $({\mathfrak a},{\mathfrak c})= 
({\mathfrak e},{\mathfrak p})$, where the empty string ${\mathfrak e}$ 
corresponds to the neutral element. This makes it possible to give simple proofs
by induction. We exemplify this in three directions: Prime ideals, upper cluster 
algebras, and the diagonal of a quantized minor.
\end{abstract}
\subjclass[2010]{MSC 17B37 (primary),\ MSC 13F60, \ MSC 16T20 (primary),  \ MSC 17A45  (secondary), \and MSC 20G42 (secondary)}
\maketitle

\section{Introduction}

We study a class of quadratic algebras connected to quantum parabolics and 
double quantum Schubert cells. We begin by considering a finite-dimensional 
simple Lie algebra ${\mathfrak g}$ over ${\mathbb C}$ and a parabolic 
sub-algebra ${\mathfrak p}\subset{\mathfrak g}$.  Then we consider a fixed Levi 
decomposition
\begin{equation}
{\mathfrak p}={\mathfrak l}+{\mathfrak u},
\end{equation}
with ${\mathfrak u}\neq 0$ and ${\mathfrak l}$  the Levi subalgebra.
  
\medskip

The main references for this study are the articles by A. Berenstein and A. 
Zelevinski \cite{bz} and by  C. Geiss, B. Leclerc, J. Schr\"oer \cite{leclerc}. 
We also refer to \cite{jak-cen} for further background.

\medskip

Let, as usual,  $W$ denote the Weyl group. Let $W_p=\{w\in W\mid 
w(\triangle^-)\cap \triangle^+\subseteq \triangle^+({\mathfrak l}) \}$ and 
$W^p$, by some called the Hasse Diagram of $G\backslash P$,  denote the usual 
set of minimal length  coset representatives of $W_p\backslash W$. Our primary 
input is a pair of Weyl group elements $w^{\mathfrak a},w^{\mathfrak c}\in W^p$ 
such that $w^{\mathfrak a}\leq w^{\mathfrak c}$. We will often, as here,  label 
our elements $w$ by ``words'' ${\mathfrak a}$; $w=w^{\mathfrak a}$,  in a 
fashion similar, though not completely identical,  to that of \cite{bz}.  
Details 
follow in later sections, but we do mention here that the element $e$ in $W$ is 
labeled by ${\mathfrak e}$ corresponding to the empty string; 
$e=\omega^{\mathfrak e}$ while the longest elements in $W^p$ is labeled by 
${\mathfrak p}$. 

To each pair $w^{\mathfrak a},w^{\mathfrak c}$ as above we construct explicitly 
a quantum seed
\begin{equation}{\mathcal Q}_q({\mathfrak a},{\mathfrak 
c}):=({\mathcal C}_q({\mathfrak a},{\mathfrak 
c}), {\mathcal L}_q({\mathfrak a},{\mathfrak 
c}), {\mathcal B}_q({\mathfrak a},{\mathfrak 
c})).\end{equation}

The cluster ${\mathcal C}_q({\mathfrak a},{\mathfrak 
c})$ generates a quadratic algebra ${\mathcal A}_q({\mathfrak a},{\mathfrak 
c})$  in the space of functions on ${\mathcal U}_q({\mathfrak n})$. 
  
\medskip

After that we define transitions 
\begin{equation}{\mathcal Q}_q({\mathfrak a},{\mathfrak 
c})\rightarrow {\mathcal Q}_q({\mathfrak a}_1,{\mathfrak 
c}_1).   
\end{equation}
We call our  transitions quantum Schur (creation/annihilation) mutations and 
prove that they are indeed just (composites of) quantum mutations in the sense 
of Berenstein and Zelevinski. These actually have to be augmented by what we 
call creation/annihilation mutations which are necessary since we have to work 
inside a larger ambient space. To keep the full generality, we may also have to 
restrict our seeds to sub-seeds.

The natural scene turns out to be   
\begin{equation}{\mathcal Q}_q({\mathfrak a},{\mathfrak b},{\mathfrak 
c}):=({\mathcal C}_q({\mathfrak a},{\mathfrak b},{\mathfrak 
c}), {\mathcal L}_q({\mathfrak a},{\mathfrak b},{\mathfrak 
c}), {\mathcal B}_q({\mathfrak a},{\mathfrak b},{\mathfrak 
c})),\end{equation}
which analogously is determined by a triple $w^{\mathfrak a},w^{\mathfrak 
b},w^{\mathfrak c}\in W^p$ such that $w^{\mathfrak a}\leq w^{\mathfrak b}\leq 
w^{\mathfrak c}$.

Later we extend our construction to even \begin{equation}
{\mathcal Q}_q({\mathfrak r}_1,\dots, {\mathfrak 
r}_{n-1},{\mathfrak r}_n)\textrm{ and }{\mathcal A}_q({\mathfrak r}_1,\dots, 
{\mathfrak 
r}_{n-1},{\mathfrak r}_n),
 \end{equation}
 though we do not use it here for anything specific.

It is a major point of this study to establish how our seeds and algebras can 
be 
constructed, inside an ambient space, starting from a single variable (indeed: 
none). In this sense the quantized generalized flag manifold of $(G/P)_q$ as 
built from quantized Schubert Cells can be built from a single cell. Furthermore, we prove that we can pass between our seeds by Schubert creation 
and 
annihilation mutations  inside a larger ambient space. 
 
 \medskip
 
 This sets the stage for (simple) inductive arguments which is a major point of this article, and is what we will 
pursue here.
 \medskip

 We first prove by induction that the two-sided sided ideal 
$I({\det}_{s}^{{\mathfrak 
a},{\mathfrak 
c}})$ in ${\mathcal A}_q({\mathfrak a},{\mathfrak c})$ generated by the quantized 
minor ${\det}_{s}^{{\mathfrak 
a},{\mathfrak 
c}}$ is  prime.

Then we prove that each upper cluster algebra ${\mathbb U}({\mathfrak 
a},{\mathfrak c})$ equals its quadratic algebra ${\mathcal A}_q({\mathfrak 
a},{\mathfrak c})$.
 
 There is a sizable overlap between these result and results previously 
obtained by K. Goodearl M. Yakimov (\cite{good},\cite{good1}).
 
\medskip 
 
 We further use our method to study the diagonal of a quantum minor.
 
 \medskip
 
The idea of induction in this context was introduced in \cite{jz} and  applications were studied in 
the case of a specific type of parabolic related to type $A_n$. Further ideas 
relating to explicit constructions of compatible pairs in special cases were 
studied in \cite{jp}.

\bigskip

\section{A little about quantum groups and cluster algebras}

\subsection{2.1 Quantum Groups}

We consider quantized enveloping algebras $U={\mathcal
U}_q({\mathfrak g})$ in
the standard notation given either eg. by
Jantzen (\cite{jan}) or by Berenstein and Zelevinsky
(\cite{bz}), though their
assumptions do not coincide completely.
To be completely on the safe side, we state our assumptions and
notation, where
it may differ: Our algebra is a Hopf
algebra defined in the usual fashion from a semi-simple
finite-dimensional
complex Lie algebra ${\mathfrak g}$. They are
algebras over ${\mathbb Q}(q)$. $\Phi$ denotes a given set of
roots and throughout, $\Pi=\{\alpha_1,\alpha_2,\dots,\alpha_R\}$ a
fixed choice of simple roots. Our
generators are then given as
$$\{E_\alpha,F_\alpha,K^\alpha\}_{\alpha\in\Pi},$$
but we will allow elements of the form  $K^\eta$ for any integer weight. $W$ denotes the Weyl group defined by $\Phi$. 

\medskip

Finally we let $\{\Lambda_\alpha\mid\alpha\in\Pi\}$ denote the set of 
fundamental
weights. We assume throughout that the diagonalizing elements $d_\alpha$ are 
determined
by
\begin{equation}
 \forall \alpha\in\Pi:(\Lambda_\alpha,\alpha)=d_\alpha.
\end{equation}

\medskip
\begin{Lem}[(2.27) in \cite{fz}]\label{3.1}  Let $\alpha_i\in \Phi$. Then 
$$(\sigma_i+1)(\Lambda_i)+\sum_{j\neq 
i}a_{ji}(\Lambda_j)=0.$$
\end{Lem}

\medskip

\subsection{Quantum Cluster Algebras}

We take over without further ado the terminology and constructions 
of
(\cite{bz}).  Results from \cite{leclerc} are also put to good use.

\begin{Def}
We say that two elements $A,B$ in some algebra over ${\mathbb C}$ $q$-commute 
if, for some $r\in{\mathbb R}$:
\begin{equation}AB=q^rBA.
\end{equation}
\end{Def}

\medskip

To distinguish between the original mutations and the more elaborate ones we need here, and to honor the founding fathers A. Berenstein, S. Fomin, 
and 
A. Zelevinski,  we use the following terminology:
\begin{Def}A quantum mutation as in \cite{bz} is called a BFZ-mutation. \end{Def}

\bigskip

\subsection{A simple observation}

If $\underline{a}=(a_1,a_2,\dots,a_{\mathfrak m})$ and
$\underline{f}=(f_1,f_2,\dots,f_{\mathfrak m})$ are vectors 
then\begin{Lem}{(\cite{jz})}\label{2.22}\begin{equation}
{\mathcal L}_q(\underline{a})^T=(\underline{f}
)^T\Leftrightarrow\forall
i:X_iX^{\underline{a}}=q^{f_i}X^{\underline{a}}X_i.\end{equation}
In particular, if there exists a $j$ such that $\forall i:
f_i=-\delta_{i,j}$
then the column vector $\underline{a}$ can be the $j$th column
in the matrix ${\mathcal B}$
of a compatible pair.
\end{Lem}
\noindent However simple this actually is, it will have a great
importance later
on.

\bigskip

\section{On Parabolics}

The 
origin of the following lies in A.
Borel \cite{borel}, and B. Kostant \cite{kos}. Other main
contributors are
\cite{bgg} and
\cite{stein}. See also \cite{cap}. We have also found (\cite{sager}) useful.

\medskip

\begin{Def} Let $w\in W$. Set $$\Phi_\omega=\{\alpha\in \Delta^+\mid
w^{-1}\alpha\in
\Delta^-\}=w( \Delta^-)\cap  \Delta^+.$$\end{Def}

\medskip

We have that $\ell(w)=\ell(w^{-1})=\vert\Phi_\omega\vert$.

\medskip

We set
$\Phi_\omega=\Delta^+(w)$. 

\medskip

From now on, we work with a fixed 
parabolic 
$\mathfrak p$ with a Levi decomposition
\begin{equation}
{\mathfrak p}={\mathfrak l}+{\mathfrak u},
\end{equation}
where ${\mathfrak l}$ is the Levi subalgebra, and where we assume ${\mathfrak 
u}\neq 0$,

\medskip

Let

\begin{Def}
\begin{eqnarray*}
W_p&=&\{w\in W\mid \Phi_\omega\subseteq \Delta^+({\mathfrak l})\},\\
W^p&=&\{w\in W\mid \Phi_\omega\subseteq \Delta^+({\mathfrak u})\}.
\end{eqnarray*}
$W^p$ is a set of distinguished representatives of the right
coset space
$W_p\backslash W$.
\end{Def}

\medskip

It is well known (see eg (\cite{sager})) that any $w\in W$ can be written 
uniquely as $w=w_pw^p$ with $w_p\in W_p$ and $w^p\in W^p$.

One defines, for each $w$ in the Weyl Group $W$, the Schubert cell $X_w$. This
is a cell in
${\mathbb P}(V)$, the projective space over a specific
finite-dimensional
representation of ${\mathfrak g}$. The closure,
${X_w}$, is called a
Schubert variety. The main classical theorems are

\begin{Thm}[Kostant,\cite{kos}]$$G/P=\sqcup_{w\in
W^p}X_w.$$\end{Thm}

{\begin{Thm}[\cite{stein}]\label{stein}
Let $w,w'\in W^p$. Then $$X_{w'}\subseteq {{X_{w}}}$$
if and only
$w'\leq w$ in the usual Bruhat ordering.
\end{Thm}

\medskip

If $\omega^{\mathfrak r}=\omega_m\tilde\omega$ and 
$\omega_{m}=\omega_n\hat\omega$ with $\omega_n,\omega_m\in W^P$
and all Weyl
group elements reduced, we say that $\omega_n<_L\omega_m$ if $\hat\omega\neq 
e$. 
This is the weak left Bruhat order.

\medskip

\section{The quadratic algebras}\label{sec4}

Let $\omega=s_{\alpha_1}s_{\alpha_2}\dots
s_{\alpha_t}$ be an element of the Weyl group written in reduced form. 
Following 
Lusztig (\cite{luz}), we construct roots $\gamma_i=\omega_{i-1}(\alpha_i)$ and 
elements $Z_{\gamma_i}\in {\mathcal U}_q({\mathfrak n}_\omega)$. 

The following result is well known, but notice a change $q\to q^{-1}$ in 
relation to 
(\cite{jak-cen}).

\begin{Thm}[\cite{lev},\cite{lev0}] \label{4.1}Suppose that  $1\leq i<j\leq t$. Then
$$Z_{i}Z_{j}=q^{-( 
\gamma_i,\gamma_j)}Z_{j}Z_{i} + {\mathcal R}_{ij},$$ where ${\mathcal R}_{ij}$ 
is of lower order in the sense that it involves only elements $Z_k$ with $i< k<
j$. 
Furthermore, the elements 
$$Z_t^{a_t}\dots Z_2^{a_2}Z_1^{a_1}$$ with $a_1,a_2,\dots,a_t\in{\mathbb N}_0$ 
form a basis of ${\mathcal U}_q({\mathfrak n}_\omega)$.
\end{Thm}
Our statement follows \cite{jan},\cite{jan2}. Other authors, eg.
\cite{lev},
\cite{leclerc} have used the other Lusztig braid operators. The
result is
just a difference between $q$ and $q^{-1}$. Proofs of this
theorem which are
more accessible are available (\cite{cp},\cite{jan2}).

\medskip

It is known that this algebra is isomorphic to the algebra of functions on 
${\mathcal U}_q({\mathfrak n}_\omega)$ satisfying the usual finiteness 
condition. It is analogously equivalent to the algebra of functions on 
${\mathcal U}^-_q({\mathfrak n}_\omega)$ satisfying a  similar finiteness 
condition. See eg (\cite{leclerc}) and (\cite{jan}).
\bigskip

\section{basic structure}\label{sec5}

Let $\omega^{\mathfrak p}$ be the maximal element in $W^p$. It
is the one which
maps all roots in $\Delta^+({\mathfrak
u})$ to $\Delta^-$. (Indeed: To $\Delta^-({\mathfrak u})$.)
Let $w_0$ be the
longest element in $W$ and $w_L$ the
longest in the Weyl group of ${\mathfrak l}$, Then
\begin{equation}w^{\mathfrak p}w_L=w_0.\end{equation}

Let 
$\omega^{\mathfrak 
r}=\sigma_{i_1}\sigma_{i_2}\cdots\sigma_{i_r}\in W^p$ be written in 
a fixed reduced form. Then
$\ell(\omega^{\mathfrak 
r})=r$. We assume here that $r\geq 1$. We set $e=\omega^{\mathfrak e}$ and 
$\ell(\omega^{\mathfrak e})=0$ where ${\mathfrak e}$ 
denotes the empty set, construed as the empty sequence. We also let ${\mathfrak r}$ denote the sequence 
$i_1,i_2,\dots,i_r$ if ${\mathfrak r}\neq {\mathfrak e}$. If a sequence ${\mathfrak s}$ corresponds to 
an analogous element $\omega^{\mathfrak s}\in W^p$ we define
\begin{equation}
{\mathfrak s}\leq {\mathfrak r}\Leftrightarrow \omega^{\mathfrak s}\leq_L 
\omega^{\mathfrak r}.
\end{equation}

 Set  
\begin{equation}\Delta^+(\omega^{\mathfrak 
r})=\{ \beta_{i_1},\dots,\beta_{i_r}\}.\end{equation}

\medskip

\begin{Def}
Let
${\mathbf b}$
denote the
map $\Pi\to\{1,2,\dots,R\}$ defined by ${\mathbf
b}(\alpha_i)=i$.
Let  $\overline\pi_{\mathfrak r}:\{1,2,\dots,
r\}\to\Pi$ be given by
\begin{equation}\overline\pi_{\mathfrak
r}(j)=\alpha_{i_j}.\end{equation}
If $\overline\pi_{\mathfrak r}(j)=\alpha$ we say that $\alpha$
(or
$\sigma_\alpha$) occurs at position 
$j$ 
in $w^{\mathfrak r}$, and we say that
$\overline\pi_{\mathfrak r}^{-1}(\alpha)$ are the positions at
which $\alpha$
occurs in $w$.
Set \begin{equation}
{\pi}_{\mathfrak r}={\mathbf b}\circ\overline\pi_{\mathfrak r}.       
\end{equation}
\end{Def}

$\pi_{\mathfrak e}$ is construed as a map whose image is the empty set.

\bigskip

Recall from (\cite{jak-cen}):

\medskip

\begin{Def}Let $\omega^{\mathfrak r}\in W^p$ be given and suppose  $s\in 
Im(\pi_{\mathfrak r})$.  Then
$s=\pi_{\mathfrak r}(n)$ for some $n$ and we set 
$\omega_n:=\sigma_{i_1}\sigma_{i_2}\cdots\sigma_{i_n}$.  Suppose
$\omega_n=\omega_1
\sigma_{i_n}\omega_2\dots\omega_t
\sigma_{i_n}$ and $\omega_i\in
W\setminus\{e\}$ for $i>1$. Further assume that each $\omega_i$ is reduced and 
does 
not contain any $\sigma_{i_n}$. We denote this simply  as $n\leftrightarrow 
(s,t)$. We further write 
$\beta_{n}\leftrightarrow \beta_{s,t}$ and 
\begin{equation}\omega_n\leftrightarrow \omega_{s,t}\end{equation} if $n,s,t$ 
are
connected as
above. It is convenient to set $\omega_{s,0}=e$ for all $s\in\{1,2,\dots, R\}$.

For a fixed $s\in\{1,2,\dots,R\}$ we let $s_{\mathfrak r}$ denote the
maximal such $t$. If there is no such decomposition we set $t=0$. So, in 
particular, $s_{\mathfrak e}=0$, and $s_{\mathfrak r}$ is the number of times $\sigma_s$ occurs in $\omega^{\mathfrak
r}$. 
Finally we set (cf. (\cite{jak-cen}))
\begin{equation}
{\mathbb U}({\mathfrak r})=\{(s,t)\in {\mathbb N}\times {\mathbb
N}_0\mid
1\leq s\leq 
R\textrm{ 
and }0\leq t\leq s_{\mathfrak r}\}.
\end{equation}
\end{Def}

\medskip

Notice that if $(s,t)\in{\mathbb U}({\mathfrak r})$ then we may
construct a subset
${\mathbb U}({\mathbf s}, {\mathbf t})$ of ${\mathbb U}$ by the above recipe,
replacing
$\omega^{\mathfrak r}$ by $\omega_{s,t}$. In this subset $t$ is
maximal. Likewise, if ${\mathfrak s}\leq {\mathfrak r}$ we have of course 
${\mathbb U}({\mathfrak s})\subseteq {\mathbb U}({\mathfrak r})$ and may set 
${\mathbb U}({\mathfrak r}\setminus {\mathfrak s})={\mathbb U}({\mathfrak 
r})\setminus {\mathbb U}({\mathfrak s})$.

\bigskip

\section{Key structures and background results}

\subsection{Quantized minors}

Following a construction of classical minors by S. Fomin and
A. Zelevinsky
\cite{fz}, the last mentioned and A. Berenstein have introduced a
family of quantized
minors $\Delta_{u\cdot\lambda,v\cdot\lambda}$ in \cite{bz}. These
are elements of
the quantized coordinate ring ${\mathcal O}_q(G)$. The results by K. Brown and 
K. Goodearl (\cite{brown}) were important in this process.

{The element 
$\Delta_{u\cdot\lambda,v\cdot\lambda}$ is determined by $u,v\in W$ and a 
positive weight $\lambda$. We will always assume that $u\leq_L v$.

\medskip

\subsection{Identifications} There is a well-known pairing
between
${\mathcal 
U}^{\leq}$ and ${\mathcal U}^{\geq}$ (\cite{jan}) and there is a
unique
bilinear 
form on ${\mathcal U}_q({\mathfrak n})$. With this we can
identify $({\mathcal
U}^{\geq})^*$ with ${\mathcal U}^{\geq}$. 
One can even define a product in $({\mathcal U}_q({\mathfrak
n}))^*$ that makes
it isomorphic to ${\mathcal U}_q({\mathfrak n})$ \cite{leclerc}.
We can in this
way identify the elements $\Delta_{u\cdot\lambda,v\cdot\lambda}$
with elements of ${\mathcal U}^{\geq}$.

\subsection{Key results from \cite{bz} and \cite{leclerc}}

The 
quantized minors are by definitions functions on ${\mathcal
U}_q({\mathfrak g})$
satisfying certain finiteness conditions. 
What is needed first are certain commutation relations that they
satisfy.
Besides this, they can be restricted to being functions on
${\mathcal
U}_q({\mathfrak b})$ and even on ${\mathcal U}_q({\mathfrak
n})$. Our main references here are (\cite{bz}) and (\cite{leclerc}); the 
details of the following can be found in the latter.

\medskip

\begin{Lem}[\cite{bz}]The
element $\triangle_{u\lambda,v\lambda}$ indeed depends only on
the weights
$u\lambda,v\lambdaλ$, not on the choices of
$u, v$ and their reduced words.
\end{Lem}

\medskip

\begin{Thm}[A version of Theorem~10.2 in \cite{bz}]
\label{10.2}For any
$\lambda,\mu\in P^+$, and $s, s', t, t' \in W$ such that
$$\ell(s's) = \ell(s') + \ell(s), \ell(t't) = \ell(t') + \ell(t)
,$$the
following holds:
$$ \triangle_{s's\lambda,t'\lambda} · \triangle_{s'\mu,t't\mu} =q^{(s\lambda | 
\mu) - (\lambda |
t\mu)}\triangle_{s'\mu,t't\mu} ·
\triangle_{s's\lambda,t'\lambda}.$$
\end{Thm}
\smallskip

It is very important for the following that the conditions essentially are on the Weyl group elements. The requirement
on $\lambda,\mu$ is furthermore independent of those.

\bigskip

An equally important fact we need is the following $q$-analogue of 
\cite[Theorem~1.17]{fz}:

\begin{Thm}[\cite{leclerc}, Proposition~3.2]\label{3.2}
Suppose that for $u,v\in W$ and $i\in I$ we have
$l(us_i)=l(u)+1$ and
$l(vs_i)=l(v)+1$.  Then

\begin{equation}\label{eq3.2}
\Delta_{us_i(\Lambda_i),vs_i(\Lambda_i)}\,\Delta_{u(\Lambda_i),v(\Lambda_i)}=
({q^{-d_i}})\Delta_{us_i(\Lambda_i),v(\Lambda_i)}\,\Delta_{
u(\Lambda_i),
vs_i(\Lambda_i)}+
\prod_{j\neq i}\Delta_{u(\Lambda_j),v(\Lambda_j)}^{-a_{ji}}
\end{equation}
holds in ${\mathcal O}_q(\frak g)$. 
\end{Thm}

(That a factor $q^{-d_i}$ must be inserted for the general case is clear.)

\medskip

One considers in \cite{leclerc}, and transformed to our terminology, modified 
elements 
\begin{equation}\label{59}D_{\xi,\eta}=\triangle_{\xi,\eta}K^{-\eta}.\end{equation}
We suppress here the restriction map $\rho$, and our  $K^{-\eta}$ is denoted as
$\triangle^\star_{\eta,\eta}$ in \cite{leclerc}. The crucial property is that
\begin{equation}
K^{-\eta}\triangle_{\xi_1,\eta_1}=q^{-(\eta,\xi_1-\eta_1)}\triangle_{\xi_1,
\eta_1}K^{-\eta}.
\end{equation}

The family $D_{\xi,\eta}$
satisfies equations analogous to those in  Theorem~\ref{10.2} subject to the 
same restrictions on the relations between the weights. 

\medskip

The following result is important:

\begin{Prop}[\cite{leclerc}]
Up to a power of $q$, the following holds:
\begin{equation}Z_{c,d}=D_{\omega^{\mathfrak 
r}_{c,d-1}(\Lambda_c),\omega^{\mathfrak 
r}_{c,d}(\Lambda_c)}.
\end{equation}
\end{Prop}

\medskip

We need a small modification of the elements $D_{\xi,\eta}$ of \cite{leclerc}:

\begin{Def}
\begin{equation}E_{\xi,\eta}:=q^{\frac14(\xi-\eta,\xi-\eta)+\frac12(\rho,
\xi-\eta)}D_{\xi,\eta}.\end{equation}
\end{Def}

It is proved in  (\cite{ki}), (\cite{re})) that $E_{\xi,\eta}$ 
is invariant under the dual
bar anti-homomorphism augmented by $q\to q^{-1}$.

Notice that this change does not affect commutators:
\begin{equation}
 D_1D_2=q^\alpha D_2D_1\Leftrightarrow E_1E_2=q^\alpha E_2E_1
\end{equation}
if $E_i=q^{x_i}D_i$ for $i=1,2$.

\medskip

\begin{Def}We say that
\begin{equation}\label{less}
E_{\xi,\eta}<E_{\xi_1,\eta_1}
\end{equation}
if $\xi=s's\lambda$, $\eta=t'\lambda$, $\xi_1=s'\mu$ and $\eta_1=t't\mu$ and 
the conditions of Theorem~\ref{10.2} are satisfied. 
\end{Def}

\medskip
The crucial equation is 
\begin{Cor}
\begin{equation}
E_{\xi,\eta}<E_{\xi_1,\eta_1}\Rightarrow 
E_{\xi,\eta}E_{\xi_1,\eta_1}=q^{(\xi-\eta,\xi_1+\eta_1)}E_{\xi_1,\eta_1}E_{\xi,
\eta}.
\end{equation}\end{Cor}

\bigskip

\subsection{Connecting with the toric frames}

\begin{Def}Suppose that $\triangle_i$, $i=1,\dots,r$ is a family of mutually
$q$-commuting elements. Let $n_1,\dots,n_r\in{\mathbb Z}$. We then set
\begin{equation}N(\prod_{i=1}^r
\triangle_i^{n_i})=q^m\prod_{i=1}^r\triangle_i^{n_i},
\end{equation}where $q^m$ is determined by the requirement that 
\begin{equation}q^{-m}\triangle_r^{n_r}\dots
\triangle_2^{n_2}\triangle_1^{n_1}= q^m
\triangle_1^{n_1}\triangle_2^{n_2}\dots \triangle_r^{n_r}.
\end{equation}
\end{Def}
It is easy to see that 
\begin{equation}
\forall \mu\in S_r: N(\prod_{i=1}^r
\triangle_{\mu(i)}^{n_{\mu(i)}})=N(\prod_{i=1}^r.
\triangle_i^{n_i})
\end{equation}

It is known through \cite{bz} that eg. the quantum minors are independent of 
the 
choices of the reduced form of $\omega^{\mathfrak r}_{\mathfrak p}$. Naturally, 
this carries over to  $\omega^{\mathfrak r}$. The quadratic algebras we have 
encountered are independent of actual choices. In the coming definition we wish 
to maintain precisely the right amount of independence. 
 
Let us now formulate Theorem~\ref{3.2} in our language while using the  
language and notation of 
toric frames from \cite{bz}. In the following Theorem we first state a formula which  uses our terminology, and then we reformulate it in the 
last two lines in terms of toric frames $M$. These frames are defined by a 
cluster made up by certain elements of the form $E_{\xi,\eta}$ to be made more 
precise later.

\begin{Thm}\label{toric}
\begin{eqnarray}
 E_{us_i\Lambda_i,vs_i\Lambda_i}&=&N\left( E_{us_i\Lambda_i,v\Lambda_i}
 E_{u\Lambda_i,vs_i\Lambda_i} E_{u\Lambda_i,v\Lambda_i}^{-1}\right)
\\\nonumber&+&N\left((\prod_{j\neq i}
E_{u(\Lambda_j),v(\Lambda_j)}^{-a_{ji}}) 
E_{u(\Lambda_i),v(\Lambda_i)}^{-1}\right)\\
&=&M(E_{us_i(\Lambda_i),v(\Lambda_i)}+E_{u(\Lambda_i),
vs_i(\Lambda_i)}-E_{u(\Lambda_i),v(\Lambda_i)})\\&+&
M(\sum_{j\neq 
i}-a_{ji}E_{u(\Lambda_j),v(\Lambda_j)}-E_{u(\Lambda_i),v(\Lambda_i)}).
\end{eqnarray}
\end{Thm}

\smallskip

\noindent{\em Proof of Theorem~\ref{toric}:} We first state a lemma whose proof 
is omitted as it is  straightforward.

\begin{Lem}
Let $\Delta_{\xi_k}$ be a family of $q$-commuting elements of weights $\xi_k$,
$k=1,\dots,r$ in the sense that for any weight $b$:\begin{equation}
\forall k=1,\dots,r: K^b\Delta_{\xi_k}=q^{(b,\xi_k)}\Delta_{\xi_k}K^b.
\end{equation}
Let $\alpha$ be defined by
\begin{equation}
\Delta_{\xi_r}\cdots\Delta_{\xi_1}=q^{-2\alpha} 
\Delta_{\xi_1}\cdots\Delta_{\xi_r}
\end{equation}Furthermore, let $b_1,\dots,b_r$ be integer weights. Then
\begin{eqnarray}
&(\Delta_{\xi_1}\Delta_{\xi_2}\cdots\Delta_{\xi_r})
K^{b_1}K^{b_2}\cdots
K^{b_r}=\\\nonumber 
&q^{\sum_{k<\ell}(b_k,\xi_\ell)}(\Delta_{\xi_1}K^{b_1})(\Delta_{\xi_2}K^{
b_2})\cdots(\Delta_{\xi_r}K^{b_r}),\textrm{ and,}\\\nonumber
&(\Delta_{\xi_r}K^{b_r})\cdots(\Delta_{\xi_1}K^{b_1})=\\\nonumber &q^{-2\alpha}
q^{(\sum_{k<\ell}-\sum_{\ell<k})(b_\ell,\xi_k)}(\Delta_{\xi_1}K^{b_1}
)\cdots(\Delta_{\xi_r}K^{b_r}),
\textrm{ so that}\\\nonumber 
&(\Delta_{\xi_1}K^{b_1})\cdots(\Delta_{\xi_r}K^{b_r})=\\\nonumber &q^{\alpha}
q^{-\frac12(\sum_{k<\ell}-\sum_{\ell<k})(b_\ell,\xi_k)}N\left(
(\Delta_{\xi_1}K^{b_1})\cdots(\Delta_{\xi_r}K^{b_r})\right).
\end{eqnarray}
Finally,
\begin{eqnarray}
&q^{-\alpha}(\Delta_{\xi_1}\Delta_{\xi_2}\cdots\Delta_{\xi_r})
K^{b_1}K^{b_2} \cdots K^{b_r}=\\\nonumber&q^{-\frac12(\sum_{\ell\neq 
k})(b_\ell,\xi_k)}N\left(
(\Delta_{\xi_1}K^{b_1})\cdots(\Delta_{\xi_r}K^{b_r})\right).
\end{eqnarray}
\end{Lem}

We apply this lemma first to the case where the elements $\xi_k$ are taken from 
the set
$\{-\textrm{sign}(a_{ki})(u\Lambda_k-v\Lambda_k)\mid a_{ki}\neq0 \}$ and where 
each
element corresponding to an $a_{ki}<0$ is taken $-a_{ki}$ times. Then 
$r=\sum_{k\neq
i}\vert a_{ji}\vert+1$. The terms considered actually commute so that here, 
$\alpha=0$.
The weights $b_k$ are chosen in the same fashion, but here
$b_k=\textrm{sign}(a_{ki})(v\Lambda_k)$. We have that 
\begin{equation}\sum_{\ell\neq
k}(b_\ell,\xi_k)=\left(\sum_{\ell}b_\ell, 
\sum_{k}\xi_k\right)-\sum_{k}(b_k,\xi_k).
\end{equation}

It follows from (\ref{3.1}) that $\sum_{\ell}b_\ell=-vs_i\lambda_i$ and
$\sum_{k}\xi_k=(us_i\Lambda_i-vs_i\lambda_i)$. Now observe that for all $k$:
$-(v\Lambda_k,(u-v)\Lambda_k)=\frac12(\xi_k,\xi_k)$. Let 
$\xi_0=(us_i-vs_i)\Lambda_i$.
The individual summands in $\sum_k(b_k,\xi_k)$ can be treated analogously. 
Keeping track of the
multiplicities and signs, it follows that
\begin{eqnarray}
q^{-\alpha}(\Delta_{\xi_1}\Delta_{\xi_2}\dots\Delta_{\xi_r})K^{b_1}K^{b_2}\dots 
K^{b_r}=\\\nonumber q^{-\frac14(\xi_0,\xi_0)+\frac14\sum_k\varepsilon_k(\xi_k,\xi_k)}N\left(
(\Delta_{\xi_1}K^{b_1})\dots(\Delta_{\xi_r}K^{b_r})\right).
\end{eqnarray}

\medskip

Let us turn to the term 
\begin{equation}
q^{-d_i}\Delta_{us_i\Lambda_i,v\Lambda_i}
\Delta_{u\Lambda_i,vs_i\Lambda_i}\Delta_{u\Lambda_i,v\Lambda_i}^{-1}K^{
-vs_i\Lambda_i}.
\end{equation}
We can of course set
$K^{-vs_i\Lambda_i}=K^{-v\Lambda_i}K^{-vs_i\Lambda_i}K^{v\Lambda_i}$. 
Furthermore, it is
known (and easy to see) that
\begin{eqnarray}
&\Delta_{u\Lambda_i,v\Lambda_i}^{-1}\Delta_{u\Lambda_i,vs_i\Lambda_i}
\Delta_{us_i\Lambda_i,v\Lambda_i}=\\\nonumber&q^{-2d_i}\Delta_{us_i\Lambda_i,
v\Lambda_i}
\Delta_{u\Lambda_i,vs_i\Lambda_i}\Delta_{u\Lambda_i,v\Lambda_i}^{-1},
\end{eqnarray}
so that $\alpha=d_i$ here. We easily get again that 
$\sum_{\ell}b_\ell=-vs_i\lambda_i$
and $\sum_{k}\xi_k=(us_i\Lambda_i-vs_i\lambda_i)$.

Let us introduce elements $\tilde
E_{\xi,\eta}=q^{\frac14(\xi-\eta,\xi-\eta)}\Delta_{\xi,\eta}K^{-\eta}$. It then 
follows
that (c.f. Theorem~\ref{3.2})

\begin{eqnarray}
\tilde E_{us_i\Lambda_i,vs_i\lambda_i}&=&N\left(\tilde 
E_{us_i\Lambda_i,v\Lambda_i}
\tilde E_{u\Lambda_i,vs_i\Lambda_i}\tilde E_{u\Lambda_i,v\Lambda_i}^{-1}\right)
\\\nonumber&+&N\left((\prod_{j\neq i}\tilde
E_{u(\Lambda_j),v(\Lambda_j)}^{-a_{ji}})\tilde 
E_{u(\Lambda_i),v(\Lambda_i)}^{-1}\right).
\end{eqnarray}

The elements $E_{\xi,\eta}$ differ from the elements $\tilde E_{\xi,\eta}$ by a 
factor
which is $q$ to an exponent which is linear in the weight $(\xi-\eta)$. Hence 
an 
equation
identical to the above holds for these elements. \qed

\bigskip

\section{Compatible pairs}

\medskip

We now construct some general families of quantum clusters and quantum seeds. 
The first, simplest, and most important, correspond to double Schubert Cells:

\medskip

Let ${\mathfrak e}\leq {\mathfrak s}<{\mathfrak t}<{\mathfrak v}\leq{\mathfrak 
p}$.  

\medskip
Set 

\begin{eqnarray*}{\mathbb U}^{d,{\mathfrak t},{\mathfrak 
v}}&:=&\{(a,j)\in {\mathbb U}({\mathfrak p})\mid a_{\mathfrak t}<j\leq 
a_{\mathfrak v}\},\\
{\mathbb U}_{R<}^{d,{\mathfrak t},{\mathfrak 
v}}&:=&\{(a,j)\in {\mathbb U}({\mathfrak p})\mid a_{\mathfrak t}<j< 
a_{\mathfrak v}\},\\
{\mathbb U}^{u,{\mathfrak s},{\mathfrak 
t}}&:=&\{(a,j)\in {\mathbb U}({\mathfrak p})\mid a_{\mathfrak s}\leq j< 
a_{\mathfrak t}\},\\
{\mathbb U}_{L<}^{u,{\mathfrak s},{\mathfrak 
t}}&:=&\{(a,j)\in {\mathbb U}({\mathfrak p})\mid a_{\mathfrak s}< j< 
a_{\mathfrak t}\}.
\end{eqnarray*}

Further, set
\begin{eqnarray} {\mathbb U}^{d,{\mathfrak t}}&=&{\mathbb U}^{d,{\mathfrak 
t},{\mathfrak 
p}},\\
  {\mathbb U}^{u,{\mathfrak t}}&=&{\mathbb U}^{d,{\mathfrak e},{\mathfrak 
t}}.
 \end{eqnarray}
\medskip

It is also convenient to define 
\begin{Def}
\begin{eqnarray}E_s(i,j)&:=&E_{\omega^{{\mathfrak 
p}}_{(s,i)}\Lambda_s,\omega^{\mathfrak p}_{(s,j)}\Lambda_s} \quad(0\leq i<j\leq 
s_{{\mathfrak p}}).
\end{eqnarray}
For $j'\geq s_{\mathfrak t}$ we set
\begin{equation}
E^d_{\mathfrak t}(s,j'):=E_s(s_{\mathfrak t},j').
\end{equation}
For $j'\leq s_{\mathfrak t}$ we set
\begin{equation}
E^u_{\mathfrak t}(s,j'):=E_s(j',s_{\mathfrak t}).
\end{equation}
Finally, we set 
\begin{eqnarray}
{\mathcal C}_q^d({\mathfrak t},{\mathfrak 
v})&=&\{E^d_{\mathfrak t}(s,j')\mid (s,j')\in 
{\mathbb U}^{d,{\mathfrak t},{\mathfrak 
v}}\},\\
{\mathcal C}_q^u({\mathfrak 
s},{\mathfrak t})&=&\{E^u_{\mathfrak t}(s,j'); (s,j')\in {\mathbb 
U}^{u,{\mathfrak 
s},{\mathfrak t}}\},\\
{\mathcal C}_q^d({\mathfrak t})&=&{\mathcal C}_q^d({\mathfrak t},{\mathfrak 
p}),\textrm{ and}\\
{\mathcal C}_q^u({\mathfrak t})&=&{\mathcal C}_q^u({\mathfrak s},{\mathfrak 
t}).
\end{eqnarray}
\end{Def}

It is clear that ${\mathcal C}_q^d({\mathfrak t},{\mathfrak 
v})\subseteq {\mathcal C}_q^d({\mathfrak t})$ for any ${\mathfrak v}>{\mathfrak 
t}$ and ${\mathcal C}_q^u({\mathfrak s},{\mathfrak 
t})\subseteq {\mathcal C}_q^u({\mathfrak t})$ for any ${\mathfrak s}<{\mathfrak 
t}$.

\begin{Lem}The elements in the set ${\mathcal C}_q^d({\mathfrak t})$ are 
$q$-commuting and the elements in the set 
${\mathcal C}_q^u({\mathfrak t})$ are 
$q$-commuting.\label{above}
\end{Lem}

The proof is omitted as it is very similar to the proof of 
Proposition~\ref{7.13} which comes later.

\smallskip

\smallskip

\begin{Def}${\mathcal A}_q^d({\mathfrak t}, {\mathfrak v})$ denotes the ${\mathbb 
C}$-algebra 
generated by ${\mathcal C}_q^d({\mathfrak t},{\mathfrak v})$ and  ${\mathcal 
A}_q^u({\mathfrak s},{\mathfrak 
t})$ denotes the ${\mathbb C}$-algebra generated by ${\mathcal C}_q^u({\mathfrak 
s},{\mathfrak t})$. Further, ${\mathcal F}_q^d({\mathfrak t},{\mathfrak v})$ and 
${\mathcal F}_q^u({\mathfrak s},{\mathfrak t})$ denote the corresponding 
skew-fields of fractions. Likewise,   ${\mathbf L}_q^d({\mathfrak t},{\mathfrak v})$ and ${\mathbf L}_q^u({\mathfrak s},{\mathfrak t})$ denote the 
respective Laurent quasi-polynomial algebras. Finally, ${\mathcal 
L}_q^d({\mathfrak t},{\mathfrak v})$ and ${\mathcal L}_q^u({\mathfrak s},{\mathfrak 
t})$ denote the symplectic forms
associated with the clusters ${\mathcal C}_q^d({\mathfrak t},{\mathfrak v})$, and 
${\mathcal C}_q^u({\mathfrak s},{\mathfrak t})$, respectively.
\end{Def}

\medskip

\begin{Def}Whenever ${\mathfrak a}<{\mathfrak b}$, we set 
\begin{equation}\forall s\in Im(\pi_{\mathfrak b}): {\det}_{s}^{{\mathfrak a},{\mathfrak b}}:=E_{\omega^{\mathfrak 
a}\Lambda_s,\omega^{\mathfrak b}\Lambda_s}. \end{equation}
\end{Def}

We conclude 
in particular that

\begin{Prop}\label{quasipol}The elements ${\det}_{s}^{{\mathfrak 
t},{\mathfrak p}}$ $q$-commute with all elements 
in the algebra ${\mathcal A}_q^d({\mathfrak t})$ and the elements 
${\det}_{s}^{{\mathfrak 
e},{\mathfrak t}}$ $q$-commute with all elements 
in the algebra ${\mathcal A}_q^d({\mathfrak t})$. 
\end{Prop}

\medskip
\begin{Def}An element $C$ in a quadratic algebra  ${\mathcal A}$ that 
$q$-commutes with all the  generating  elements is said to be covariant.
\end{Def}

\medskip

As a small aside, we mention the following easy generalization of the result in 
(\cite{jak-cen}): 

\begin{Prop}It ${\mathfrak a}<{\mathfrak b}$, then
the spaces ${\mathcal A}_q^u({\mathfrak a},{\mathfrak b})$ and ${\mathcal 
A}_q^d({\mathfrak a},{\mathfrak b})$ are quadratic algebras. In both cases, the 
center is given by $Ker(\omega^{\mathfrak a}+\omega^{\mathfrak b})$. The semi-group of covariant elements in generated by $\{ {\det}_{s}^{{\mathfrak a},{\mathfrak b}}\mid s\in Im(\pi_{\mathfrak b})\}$.
\end{Prop}
\medskip

We now construct some elements in ${\mathbf L}_q^d({\mathfrak t})$ and ${\mathbf 
L}_q^u({\mathfrak t})$ of fundamental importance. They are indeed monomials in 
the elements of $\left[{\mathcal C}_q^d({\mathfrak t})\right]^{\pm1}$ and   
$\left[{\mathcal C}_q^u({\mathfrak t})\right]^{\pm1}$, respectively.

\medskip

First a technical definition:

\begin{Def}
$p(a,j,k)$ denotes the largest non-negative integer for which 
$$\omega^{\mathfrak p}_{(k,p(a,j,k))}\Lambda_k=\omega^{\mathfrak 
p}_{(a,j)}\Lambda_k.$$
\end{Def}

\medskip

We also allow $E_a(j,j)$ which is defined to be $1$.

\medskip

Here are then the first building blocks:

\begin{Def}
\begin{eqnarray}\nonumber
&\forall (a,j)\in {\mathbb U}^{d,{\mathfrak t}}:\\& H^d_{\mathfrak 
t}(a,j):=E_a(a_{\mathfrak t},j)E_a(a_{\mathfrak t},j-1)
\prod_{a_{ka}<0}E_k(k_{\mathfrak t},p(a,j,k))^{a_{ka}} \\\nonumber
&\forall (a,j)\in  {\mathbb U}^{d,{\mathfrak t}}\textrm{ with }j<a_{\mathfrak 
p}:\\ 
&B^d_{\mathfrak t}(a,j):=H^d_{\mathfrak t}(a,j)(H^d_{\mathfrak t}(a,j+1))^{-1}.
\end{eqnarray}
\end{Def}
The terms $E(k_{\mathfrak t},p(a,j,k))$ and $E_a(a_{\mathfrak t},j-1)$ are 
well-defined but may become equal to $1$. Also notice that, where defined,   
$H^d_{\mathfrak t}(a,j), B^d_{\mathfrak t}(a,j)\in {\mathbf  L}_q^d({\mathfrak 
t})$.
\bigskip

\begin{Lem}\label{7.10}If $E_{\xi,\eta}<H^d_{\mathfrak t}(a,j)$ in the sense 
that it is less 
than or equal to each factor $E_{\xi_1,\eta_1}$ of $H^d_{\mathfrak t}(a,j)$ 
(and 
$<$  is defined in (\ref{less})), then 
\begin{equation}\label{54}
E_{\xi,\eta}H^d_{\mathfrak t}(a,j)=q^{(\xi-\eta,\omega^{\mathfrak 
t}(\alpha_a))}H^d_{\mathfrak t}(a,j)E_{\xi,\eta}.
\end{equation}
If $E_{\xi,\eta}\geq H^d_{\mathfrak t}(a,j)$, then 
\begin{equation}\label{55}
E_{\xi,\eta}H^d_{\mathfrak t}(a,j)=q^{(-\xi-\eta,\omega^{\mathfrak 
t}(\alpha_a))}H^d_{\mathfrak t}(a,j)E_{\xi,\eta}.
\end{equation}
\end{Lem}

\proof This follows from (\ref{less}) by observing that we have the following 
pairs $(\xi_1,\eta_1)$ occurring in $H^d_{\mathfrak t}(a,j)$: 
$$(\omega^{\mathfrak t}\Lambda_a,\omega(a,j)\Lambda_a),(\omega^{\mathfrak 
t}\Lambda_a,\omega(a,j)\sigma_a\Lambda_a),$$ and $$(-\omega^{\mathfrak 
t}\Lambda_k,-\omega(a,j)\Lambda_k) \textrm{ with multiplicity }(-a_{ka}).$$
Furthermore, as in (\ref{3.1}), 
$\Lambda_a+\sigma_a\Lambda_a+\sum_ka_{ka}\Lambda_k=0$ and, 
equivalently, $2\Lambda_a+\sum_ka_{ka}\Lambda_k=\alpha_a$ . \qed

\bigskip

\begin{Prop}\label{7.10}$\forall (a,j),(b,j')\in {\mathbb U}^{d,{\mathfrak t}}, 
j<a_{\mathfrak p}$ the following holds:
\begin{equation}
E^d_{\mathfrak t}(b,j')B^d_{\mathfrak 
t}(a,j)=q^{-2(\Lambda_a,\alpha_a)\delta_{j,j'}\delta_{a,b}}B^d_{\mathfrak 
t}(a,j)E^d_{\mathfrak t}(b,j').
\end{equation}
\end{Prop}

\medskip

\proof It is clear from the formulas (\ref{54}-\ref{55}) that if an element 
$E_{\xi,\eta}$ either is bigger than all factors in $B^d_{\mathfrak 
t}(s,j)$ or smaller than all factors, then it commutes with this element. The 
important fact now is that the ordering is independent of the fundamental 
weights $\Lambda_i$ - it depends only on the Weyl group elements. The factors 
in any $H_{\mathfrak t}^d$ are, with a fixed ${\mathfrak t}$, of the form 
$E_{\omega^{\mathfrak t}\Lambda_i,\omega\Lambda_i}$ or $E_{\omega^{\mathfrak 
t}\Lambda_a,\omega\circ\sigma_a\Lambda_a}$ for some $\omega\geq 
\omega^{\mathfrak t}$. The elements $E_{\xi,\eta}=E^d_{\mathfrak t}(b,j')$ we 
consider thus satisfy the first or the second case in Lemma~\ref{7.10} 
for either terms $H^d_{\mathfrak t}(a,j)$ and $H^d_{\mathfrak t}(a,j+1)$. 
Clearly, we then need only consider the  in-between case $H^d_{\mathfrak t}(a,j)\leq E_{\xi,\eta}\leq H^d_{\mathfrak t}(a,j+1)$, and here there appears a 
factor $q^{-2(\xi,\omega^{\mathfrak t}(\alpha_a))}$ in the commutator with 
$\xi=\omega^{\mathfrak t}\Lambda_{b}$. This accounts for the term 
$-2(\Lambda_a,\alpha_a)\delta_{a,b}$. Finally,  if $a=b$ the previous 
assumption forces $j=j'$.  \qed

\bigskip

Let us choose an enumeration
\begin{equation}
{\mathcal C}_q^d({\mathfrak t})=\{c_1,c_2,\dots, c_N\}
\end{equation}
so that each $(a,j)\leftrightarrow k$ and let us use the same enumeration of  
 the elements $B^d_{\mathfrak 
t}(a,j)$.  Set, for now $B^d_{\mathfrak 
t}(a,j)=b_k$ if $(a,j)\leftrightarrow k$. Let us also agree that the, say $n$, 
non-mutable elements $\det_s^{{\mathfrak t},{\mathfrak p}}$ of ${\mathcal 
C}_q^d({\mathfrak t})$ are written last, say numbers $N-n+1, N-n+2,\dots, N-n$. 
Then, as defined,  
\begin{equation}
\forall j=1,\dots, N-n: b_j=q^{\alpha_j}\prod_kc_k^{b_{kj}}
\end{equation}
for some integers $b_{kj}$ and some, here inconsequential, factor 
$q^{\alpha_j}$. The symplectic form yields a matrix which we, abusing notation 
slightly, also denote ${\mathcal L}_q^d({\mathfrak t})$ such that
\begin{equation}
\forall i,j=1,\dots, N:\ \left({\mathcal L}_q^d({\mathfrak 
t})\right)_{ij}=\lambda_{ij}
\end{equation}
and
\begin{equation}
\forall i,j=1,\dots, N: c_jc_j=q^{\lambda_{ij}}c_jc_i.\end{equation}
Similarly, we let ${\mathcal B}_q^d({\mathfrak t})$ denote the matrix
\begin{equation}
\forall i=1,\dots, N\ \forall j=1,\dots, N-n: \left({\mathcal B}_q^d({\mathfrak 
t})\right)_{ij}=b_{ij}.
\end{equation}
Then, where defined,
\begin{equation}
c_ib_j=\prod_kq^{\lambda_{ik}b_{kj}}b_jc_i,
\end{equation}
and Proposition~\ref{7.10} may then be restated as
\begin{equation}
\forall i=1,\dots, N\ \forall j=1,\dots, N-n: 
\sum_k\lambda_{ik}b_{kj}=-2(\Lambda_s,\alpha_s)\delta_{ij},
\end{equation}
where we assume that $i\leftrightarrow (s,\ell)$.

\medskip

We have then established
\begin{Thm}
The pair $({\mathcal L}_q^d({\mathfrak t}),{\mathcal B}_q^d({\mathfrak t}))$ is a 
compatible pair and hence,
\begin{equation}
{\mathcal Q}_q^d({\mathfrak t}):=({\mathcal C}_q^d({\mathfrak t}),{\mathcal 
L}_q^d({\mathfrak t}),{\mathcal B}_q^d({\mathfrak t}))
\end{equation}
is a quantum seed with the $n$ non-mutable elements $\det_s^{{\mathfrak 
t},{\mathfrak p}}$, $(s,s_{\mathfrak p})\in {\mathbb U}^d({\mathfrak t})$. The 
entries of the diagonal of the matrix $\tilde D=({\mathcal B}_q^d({\mathfrak 
t}))^T{\mathcal L}_q^d({\mathfrak t})$ are in the set 
$\{2(\Lambda_s,\alpha_s)\mid s=1,\dots,R\}$.
\end{Thm}

\bigskip

It ${\mathfrak v}>{\mathfrak t}$, we let  $({\mathcal L}_q^d({\mathfrak 
t},{\mathfrak v}),{\mathcal B}_q^d({\mathfrak t},{\mathfrak v}))$ denote the part 
of the compatible pair  $({\mathcal L}_q^d({\mathfrak t}),{\mathcal 
B}_q^d({\mathfrak t}))$ that corresponds to the cluster ${\mathcal 
C}_q^d({\mathfrak t},{\mathfrak v})$ and we let ${\mathcal Q}_q^d({\mathfrak 
t},{\mathfrak v})$ be the corresponding triple. It is then obvious by simple 
restriction,  that we in fact have obtained

\begin{Thm}
The pair $({\mathcal L}_q^d({\mathfrak t},{\mathfrak v}),{\mathcal 
B}_q^d({\mathfrak t},{\mathfrak v}))$ is a compatible pair and hence,
\begin{equation}
{\mathcal Q}_q^d({\mathfrak t},{\mathfrak v}):=({\mathcal C}_q^d({\mathfrak 
t},{\mathfrak v}),{\mathcal L}_q^d({\mathfrak t},{\mathfrak v}),{\mathcal 
B}_q^d({\mathfrak t},{\mathfrak v}))
\end{equation}
is a quantum seed with the $n$ non-mutable elements $\det_s^{{\mathfrak 
t},{\mathfrak v}}$, $(s,s_{\mathfrak v})\in {\mathbb U}^d({\mathfrak 
t},{\mathfrak v})$. 
\end{Thm}

\bigskip

The case of ${\mathcal C}_q^u({\mathfrak t})$ is completely analogous: Define

\begin{Def}
\begin{eqnarray}
H^u_{\mathfrak t}(a,j)&:=&E_a(j, a_{\mathfrak t})E_a(j-1,a_{\mathfrak t})
\prod_{a_{ka}<0}E_k(p(a,j,k), k_{\mathfrak t})^{a_{ka}} \ (1\leq j<a_{\mathfrak 
t}),\nonumber\\
B^u_{\mathfrak t}(a,j)&:=&H^u_{\mathfrak t}(a,j+1)(H^u_{\mathfrak t}(a,j))^{-1} 
\ (1\leq j<a_{\mathfrak t}).
\end{eqnarray}
The terms $E(p(a,j,k),k_{\mathfrak t})$ are well-defined but may become equal 
to 
$1$. Notice also the exponents on the terms $H^u_{\mathfrak t}$.
\end{Def}

The terms $E(p(a,j,k),k_{\mathfrak t})$ are well-defined but may become equal 
to 
$1$. As defined, $H^u_{\mathfrak t}(a,j)$, and $B^u_{\mathfrak t}(a,j)$ are in 
${\mathbf L}_q^u({\mathfrak t})$.

\smallskip

\begin{Prop}$\forall (a,j),(b,j')\in {\mathbb U}^{u,{\mathfrak t}}, 1\leq j$ 
the 
following holds:

\begin{equation}
E^u_{\mathfrak t}(b,j')B^u_{\mathfrak 
t}(a,j)=q^{2(\Lambda_a,\alpha_a)\delta_{j,j'}\delta_{a,b}}B^u_{\mathfrak 
t}(a,j)E^u_{\mathfrak 
t}(b,j').
\end{equation}
\end{Prop}

\medskip

We then get in a similar way

\begin{Thm}
The pair $({\mathcal L}_q^u({\mathfrak t}),{\mathcal B}_q^u({\mathfrak t}))$ is a 
compatible pair and hence,
\begin{equation}
{\mathcal Q}_q^u({\mathfrak t}):=({\mathcal C}_q^u({\mathfrak t}),{\mathcal 
L}_q^u({\mathfrak t}),{\mathcal B}_q^u({\mathfrak t}))
\end{equation}
is a quantum seed with the $n$ non-mutable elements $\det_s^{{\mathfrak 
e},{\mathfrak t}}$, $(s,s_{\mathfrak t})\in {\mathbb U}^u({\mathfrak t})$. 
\end{Thm}

\bigskip

Naturally, we even have

\begin{Thm}
The pair $({\mathcal L}_q^u({\mathfrak s},{\mathfrak t}),{\mathcal 
B}_q^u({\mathfrak s},{\mathfrak t}))$ is a compatible pair and hence,
\begin{equation}
{\mathcal Q}_q^u({\mathfrak s},{\mathfrak t}):=({\mathcal C}_q^u({\mathfrak 
s},{\mathfrak v}),{\mathcal L}_q^u({\mathfrak s},{\mathfrak t}),{\mathcal 
B}_q^u({\mathfrak s},{\mathfrak t}))
\end{equation}
is a quantum seed with the $n$ non-mutable elements $\det_s^{{\mathfrak 
s},{\mathfrak t}}$, $(s,s_{\mathfrak s})\in {\mathbb U}^u({\mathfrak 
s},{\mathfrak t})$.   
\end{Thm}

\bigskip

We now wish to consider more elaborate seeds. The first generalization is the 
most important:

Let 
\begin{equation}
{\mathfrak e}\leq {\mathfrak a}\leq {\mathfrak 
b}\leq{\mathfrak c}\leq {\mathfrak p}, \textrm{ but } {\mathfrak a}\neq  
{\mathfrak c}.\end{equation}

\begin{eqnarray}\label{l1}
{\mathcal C}_q^d({\mathfrak a},{\mathfrak b},{\mathfrak 
c})&:=&\{E^d_{\mathfrak a}(s,j)\mid (a,j)\in ({\mathbb U}^{d,{\mathfrak 
b}}\setminus {\mathbb U}^{d,{\mathfrak 
c}})= {\mathbb U}^{d,{\mathfrak 
b},{\mathfrak c}}\},\\\label{l2}{\mathcal C}_q^u({\mathfrak 
a},{\mathfrak b},{\mathfrak 
c})&:=&\{E^u_{\mathfrak c}(s,j)\mid (s,j)\in ({\mathbb U}^{u,{\mathfrak 
b}}\setminus {\mathbb U}^{u,{\mathfrak 
a}})= {\mathbb U}^{u,{\mathfrak 
a},{\mathfrak b}}\}.
\end{eqnarray}
In (\ref{l1}), ${\mathfrak a}={\mathfrak b}$ is allowed, and in (\ref{l2}), 
${\mathfrak b}={\mathfrak c}$ is allowed.

\medskip

\begin{Def}
\begin{eqnarray*}{\mathcal C}_q({\mathfrak a},{\mathfrak b},{\mathfrak 
c})&:=&{\mathcal C}_q^d({\mathfrak a},{\mathfrak b},{\mathfrak 
c})\cup{\mathcal C}_q^u({\mathfrak a},{\mathfrak b},{\mathfrak 
b}),\\
{\mathcal C}_q^o({\mathfrak a},{\mathfrak b},{\mathfrak 
c})&:=&{\mathcal C}_q^u({\mathfrak a},{\mathfrak b},{\mathfrak 
c})\cup{\mathcal C}_q^d({\mathfrak b},{\mathfrak b},{\mathfrak 
c}).
\end{eqnarray*}
\end{Def}

\begin{Prop}\label{7.13}
The elements of ${\mathcal C}_q({\mathfrak a},{\mathfrak b},{\mathfrak 
c})$ and ${\mathcal C}_q^o({\mathfrak a},{\mathfrak b},{\mathfrak 
c})$, respectively,  $q$-commute. 
\end{Prop}

\proof The two cases are very similar, so we only prove it for the first case. 
We examine 3 cases, while using the following mild version of Theorem~\ref{10.2}: 
$
\triangle_{s's\lambda,t'\lambda}$ and $\triangle_{s'\mu,t't\mu}$ $q$ commute for any $\lambda,\mu\in P^+$, and $s, s', t, t' \in W$ for which $\ell(s's) = 
\ell(s') + \ell(s), \ell(t't) = \ell(t') + \ell(t)$.

{\bf Case 1:} $E_{\mathfrak a}^d(s,t)$ and  $E_{\mathfrak a}^d(s_1,t_1)$ for $(s,t)\in {\mathbb U}^{d, {\mathfrak b},{\mathfrak c}}$ and
$(s,t)<(s_1,t_1)$:
Set $\lambda=\Lambda_s,\mu=\Lambda_{s_1}$, $s=1,s'=\omega^{\mathfrak a}$, and
$t'=\omega^{\mathfrak c}(s,t),
t't=\omega^{\mathfrak c}(s_1,t_1)$. 

{\bf Case 2:}  $E_{\mathfrak b}^u(s,t)$ and  $E_{\mathfrak b}^u(s_1,t_1)$ for $(s,t)\in {\mathbb U}^{u, {\mathfrak a},{\mathfrak b}}$ and
$(s,t)>(s_1,t_1)$:
Set $\lambda=\Lambda_s,\mu=\Lambda_{s_1}$, $t=1$,
$t'=\omega^{\mathfrak b}$ and
$s'=\omega^{\mathfrak p}(s_1,t_1), s's=\omega^{\mathfrak r}(s,t)$.

{\bf Case 3:}  $E_{\mathfrak b}^u(s,t)$ and  $E_{\mathfrak a}^d(s_1,t_1)$ for $(s,t)\in  {\mathbb U}^{u, {\mathfrak a},{\mathfrak b}}$ and
$(s_1,t_1)\in
 {\mathbb U}^{d, {\mathfrak b},{\mathfrak c}}$: 
Set $\lambda=\Lambda_s,\mu=\Lambda_{s_1}$, $s'=\omega^{{\mathfrak a}}$,
$s=\omega^{\mathfrak p}(s,t)$,
$t'=\omega^{\mathfrak b}$ and 
$t't=\omega^{\mathfrak p}(s_1,t_1)$.
\qed

\smallskip

Notice that the ordering in ${\mathbb U}^{u, {\mathfrak a},{\mathfrak b}}$ (Case 2) is the 
opposite of that of the two other cases.

\bigskip

We also define, for ${\mathfrak a}<{\mathfrak b}$,

\begin{eqnarray}
{\mathcal C}_q^u({\mathfrak a},{\mathfrak b})&=&{\mathcal C}_q^u({\mathfrak a}, 
{\mathfrak b},{\mathfrak b}),\textrm{ and}\\\nonumber {\mathcal C}_q^d({\mathfrak 
a},{\mathfrak b})&=&{\mathcal C}_q^d({\mathfrak a}, {\mathfrak a},{\mathfrak b}).
\end{eqnarray}

\medskip

We let ${\mathcal L}_q({\mathfrak a},{\mathfrak b},{\mathfrak 
c})$ and ${\mathcal L}_q^o({\mathfrak a},{\mathfrak b},{\mathfrak 
c})$ denote the corresponding symplectic matrices. We proceed to construct 
compatible pairs and give the details for just ${\mathcal C}_q({\mathfrak 
a},{\mathfrak b},{\mathfrak 
c})$.  We will be completely explicit except in the special cases 
$E^u_{\omega^{{\mathfrak 
a}}\Lambda_s,\omega^{{\mathfrak b}}\Lambda_s}={\det}_{s}^{{\mathfrak 
a},{\mathfrak b}}$ where we only give a recipe for $ {B}_q^{{\mathfrak 
a},{\mathfrak b},{\mathfrak 
c}}(s,s_{{\mathfrak a}}) $. Notice, however, the remark following (\ref{77}).

\medskip

\begin{equation}\label{72}{B}_q^{{\mathfrak a}, {\mathfrak 
b},{\mathfrak c}}(s,j):=\left\{\begin{array}{lll}B^d_{{\mathfrak 
a}}(s,j)&\textrm{if }(s,j)\in {\mathbb U}_{R<}^{d, {\mathfrak b},{\mathfrak 
c}}\\ \ \\
B^u_{{\mathfrak b}}(s,j)&\textrm{if }(s,j)\in {\mathbb U}_{L<}^{u, {\mathfrak 
a},{\mathfrak b}}  \end{array}\right..\end{equation}

\smallskip

We easily get from the preceding propositions:

\begin{Prop}\label{7.20}
Let $E(b,j')\in {\mathcal C}_q({{\mathfrak a}, {\mathfrak 
b},{\mathfrak c}})$ and let ${B}_q^{{\mathfrak a}, {\mathfrak 
b},{\mathfrak c}}(s,j)$ be as in the previous equation. Then
\begin{equation}
E(b,j'){B}_q^{{\mathfrak a}, {\mathfrak 
b},{\mathfrak 
c}}(s,j)=q^{{-2(\Lambda_s,\alpha_s)\delta_{j,j'}\delta_{s,b}}}{B}_q^{{\mathfrak 
a}, {\mathfrak 
b},{\mathfrak c}}(s,j)E(b,j'),
\end{equation}
and ${B}_q^{{\mathfrak a}, {\mathfrak 
b},{\mathfrak c}}(s,j)$ is in the algebra ${\mathcal A}_q({{\mathfrak a}, 
{\mathfrak 
b},{\mathfrak c}})$ generated by the elements of ${\mathcal C}_q({{\mathfrak a}, 
{\mathfrak 
b},{\mathfrak c}})$.
\end{Prop}

This then leaves the positions $(s,s_{\mathfrak c})\in {\mathbb 
U}^{d,{\mathfrak 
b},{\mathfrak c}}$ and $(s,s_{\mathfrak a})\in {\mathbb U}^{u,{\mathfrak 
a},{\mathfrak b}}$ to be considered. Here, the first ones are considered as the non-mutable elements. In the ambient space  ${\mathcal A}_q({{\mathfrak a}, 
{\mathfrak b},{\mathfrak c}})$, the positions in remaining cases define elements that are, in general, mutable.   

The elements in these cases are of the form $E_{\omega^{{\mathfrak 
a}}\Lambda_s,\omega^{{\mathfrak b}}\Lambda_s}$ for some $s$. To give a recipe  we define the following elements in ${\mathcal 
A}_q({{\mathfrak a}, {\mathfrak 
b},{\mathfrak c}})$:
\begin{eqnarray}&\tilde {B}_q^{{\mathfrak a},{\mathfrak b},{\mathfrak 
c}}(s,s_{{\mathfrak a}})  :=\\&\left(H^u_{{\mathfrak 
b}}(s,s_{{\mathfrak a}}+1) H^d_{{\mathfrak a}}(s,s_{{\mathfrak 
b}}+1)\right)^{-1} E_s(s_{{\mathfrak a}},s_{{\mathfrak 
b}})^2\prod_{a_{ka}<0}E_k(k_{{\mathfrak a}},k_{{\mathfrak 
b}})^{a_{ks}}.\nonumber\end{eqnarray}

If $\omega(s,s_{{\mathfrak a}}+1)=\omega^{{\mathfrak a}}\omega_x\sigma_s$ and 
$\omega(s,s_{{\mathfrak b}}+1)=\omega^{{\mathfrak b}}\omega_y\sigma_s$, and if we set $u=\omega^{{\mathfrak a}}\omega_x$, $v=\omega^{{\mathfrak b}}\omega_y$ this takes the simpler form  

\begin{equation}\label{75}\tilde {B}_q^{{\mathfrak a},{\mathfrak b},{\mathfrak 
c}}(s,s_{{\mathfrak a}}) =E_{u\sigma_s\Lambda_s,v\Lambda_s}^{-1}E_{
u\Lambda_s,v\sigma_s\Lambda_s}^{-1}\prod_{a_{ks<0}} 
E_{\omega^{\mathfrak a}\Lambda_k,v\Lambda_k}^{-a_{ks}}\prod_{a_{ks<0}}  
E_{u\Lambda_k,\omega^{\mathfrak b}\Lambda_k}^{-a_{ks}}\prod_{a_{ks<0}}  
E_{\omega^{\mathfrak a}\Lambda_k,\omega^{\mathfrak 
b}\Lambda_k}^{a_{ks}}.\end{equation}

\begin{Prop}
\begin{equation}
\forall\ell: E_{\omega^{\mathfrak a}\Lambda_\ell, \omega^{\mathfrak 
b}\Lambda_\ell}\tilde {B}_q^{{\mathfrak a},{\mathfrak b},{\mathfrak 
c}}(s,s_{{\mathfrak a}})=q^{-2\delta_{\ell.s}(\lambda_s,\alpha_s)} \tilde 
{B}_q^{{\mathfrak a},{\mathfrak b},{\mathfrak 
c}}(s,s_{{\mathfrak a}})E_{\omega^{\mathfrak a}\Lambda_\ell, \omega^{\mathfrak 
b}\Lambda_\ell}.
\end{equation}
Besides this,  $\tilde {B}_q^{{\mathfrak a},{\mathfrak b},{\mathfrak 
c}}(s,s_{{\mathfrak a}})$ commutes with everything in the cluster except possibly elements of 
the form 
$$ E_{\omega^{{\mathfrak a}}\Lambda_\ell,\omega^{{\mathfrak 
b}}\tilde\omega_y\Lambda_\ell}, \textrm{ and } E_{\omega^{{\mathfrak 
a}}\tilde\omega_x\Lambda_\ell,\omega^{{\mathfrak 
b}}\Lambda_\ell}, $$ with $1<\tilde\omega_x<\omega_x$ and     
$1<\tilde\omega_y<\omega_y$. 
\end{Prop}

The exceptional terms above are covered by Proposition~\ref{7.20} which means 
that we can in principle make a modification $\tilde {B}_q^{{\mathfrak 
a},{\mathfrak b},{\mathfrak 
c}}(s,s_{{\mathfrak a}})\to {B}_q^{{\mathfrak a},{\mathfrak b},{\mathfrak 
c}}(s,s_{{\mathfrak a}})$ where the latter expression commutes with everything 
except  $E_{\omega^{{\mathfrak a}}\Lambda_s,\omega^{{\mathfrak 
b}}\Lambda_s}$ where we get a factor $q^{-2(\Lambda_s,\alpha_s)}$.

\medskip

If $\omega_y=1$ we get a further 
simplification where now $u=\omega^{{\mathfrak a}}\omega_x$ and 
$v=\omega^{{\mathfrak 
b}}$:
\begin{equation}\label{77}\tilde {B}_q^{{\mathfrak a},{\mathfrak b},{\mathfrak 
c}}(s,s_{{\mathfrak a}}) =E_{u\sigma_s\Lambda_s,v\Lambda_s}^{-1}E_{
u\Lambda_s,v\sigma_s\Lambda_s}^{-1}\prod_{a_{ks<0}} 
E_{u\Lambda_k,v\Lambda_k}^{-a_{ks}}.\end{equation}

Here we actually have $\tilde {B}_q^{{\mathfrak a},{\mathfrak b},{\mathfrak 
c}}(s,s_{{\mathfrak a}}) ={B}_q^{{\mathfrak a},{\mathfrak b},{\mathfrak 
c}}(s,s_{{\mathfrak a}})$,  and this expression has the exact form needed for the purposes of the next section.

\medskip

We let ${\mathcal B}_q({\mathfrak a},{\mathfrak b},{\mathfrak 
c})$ and ${\mathcal B}_q^o({\mathfrak a},{\mathfrak b},{\mathfrak 
c})$ denote the corresponding symplectic matrices and can now finally define 
our quantum seeds:

\begin{equation}{\mathcal Q}_q({\mathfrak a},{\mathfrak b},{\mathfrak 
c}):=({\mathcal C}_q({\mathfrak a},{\mathfrak b},{\mathfrak 
c}), {\mathcal L}_q({\mathfrak a},{\mathfrak b},{\mathfrak 
c}), {\mathcal B}_q({\mathfrak a},{\mathfrak b},{\mathfrak 
c})).\end{equation}

\begin{Def}\begin{equation}{\mathcal Q}_q^o({\mathfrak a},{\mathfrak 
b},{\mathfrak 
c}):=({\mathcal C}_q^o({\mathfrak a},{\mathfrak b},{\mathfrak 
c}), {\mathcal L}_q^o({\mathfrak a},{\mathfrak b},{\mathfrak 
c}), {\mathcal B}_q^o({\mathfrak a},{\mathfrak b},{\mathfrak 
c})).\end{equation}\end{Def}

\medskip

According to our analysis above we have established

\begin{Thm}\label{seedth}
They are indeed seeds. The non-mutable elements are in both cases the elements 
${\det}_{s}^{{\mathfrak a},{\mathfrak c}}; s\in Im(\pi_{\omega^{\mathfrak c}})$.
\end{Thm}

\medskip

Let us  finally consider a  general situation where we are given a finite 
sequence of elements $\{\omega^{{\mathfrak r}_i}\}_{i=1}^n\in W^p$ such that 
\begin{equation}\label{genseq}
{\mathfrak e}\leq {{\mathfrak r}_1}<\dots<{{\mathfrak 
r}_n}\leq {{\mathfrak p}}.
\end{equation}
Observe that \begin{equation}\forall(s,t)\in{\mathbb U}({\mathfrak 
r}_k):\omega^{{\mathfrak r}_k}_{(s,t)}=\omega^{{\mathfrak p}}_{(s,t)}.
\end{equation}

It may of course well happen that for some $a$, and some $ {{\mathfrak 
r}_i}<{{\mathfrak 
r}_j}$,
\begin{equation}\omega^{{\mathfrak r}_i}\Lambda_a=\omega^{{\mathfrak 
r}_j}\Lambda_a.
\end{equation}

\begin{Def}Given (\ref{genseq}) we define

\begin{eqnarray}{\mathcal C}_q({\mathfrak r}_1,\dots, {\mathfrak 
r}_{n-1},{\mathfrak r}_n)&=&{\mathcal C}_q^d({\mathfrak r}_1, {\mathfrak 
r}_{n-1},{\mathfrak r}_n)\cup {\mathcal C}_q^u({\mathfrak r}_1,{\mathfrak 
r}_{2},{\mathfrak r}_{n-1})\cup\dots\\&=&\bigcup_{0<2i\leq n}{\mathcal 
C}_q^d({\mathfrak r}_i,{\mathfrak r}_{n-i},{\mathfrak r}_{n-i+1})\cup 
\bigcup_{0<2j\leq n-1}{\mathcal C}_q^u({\mathfrak r}_j,{\mathfrak 
r}_{j+1},{\mathfrak r}_{n-j}). \nonumber
\end{eqnarray}

It is also convenient to consider

\begin{eqnarray}{\mathcal C}_q^o({\mathfrak r}_1,\dots, {\mathfrak 
r}_{n-1},{\mathfrak r}_n)&=&{\mathcal C}_q^u({\mathfrak r}_1, {\mathfrak 
r}_{2},{\mathfrak r}_n)\cup {\mathcal C}_q^d({\mathfrak r}_2,{\mathfrak 
r}_{n-1},{\mathfrak r}_{n})\cup\dots\\&=&\bigcup_{0<2i\leq n}{\mathcal 
C}_q^u({\mathfrak r}_i,{\mathfrak r}_{i+1},{\mathfrak r}_{n-i+1})\cup 
\bigcup_{0<2j\leq n-1}{\mathcal C}_q^d({\mathfrak r}_{j+1},{\mathfrak 
r}_{n-j},{\mathfrak r}_{n-j+1}). \nonumber
\end{eqnarray}
\end{Def}

Notice that 
\begin{eqnarray}{\mathcal C}_q({\mathfrak r}_1,\dots, {\mathfrak 
r}_{n-1},{\mathfrak r}_n)&=&{\mathcal C}_q^d({\mathfrak r}_1, {\mathfrak 
r}_{n-1},{\mathfrak r}_n)\cup {\mathcal C}_q^o({\mathfrak r}_1,\dots, {\mathfrak 
r}_{n-2},{\mathfrak r}_{n-1})\\\nonumber
{\mathcal C}_q^o({\mathfrak r}_1,\dots, {\mathfrak 
r}_{n-1},{\mathfrak r}_n)&=&{\mathcal C}_q^u({\mathfrak r}_1, {\mathfrak 
r}_{2},{\mathfrak r}_n)\cup {\mathcal C}_q({\mathfrak r}_2,\dots, {\mathfrak 
r}_{n-1},{\mathfrak r}_n)
\end{eqnarray}

\medskip

For  the last equations, notice that ${\mathcal C}_q^{d}({\mathfrak 
e},{\mathfrak r},{\mathfrak r}) =\emptyset={\mathcal C}_q^{d}({\mathfrak 
r},{\mathfrak r},{\mathfrak r})$.

\bigskip

\begin{Prop}\label{qcom}
The spaces
\begin{equation}
{\mathcal C}_q^o({\mathfrak r}_1,\dots, {\mathfrak 
r}_{n-1},{\mathfrak r}_n)\textrm{ and }{\mathcal C}_q({\mathfrak r}_1,\dots, 
{\mathfrak 
r}_{n-1},{\mathfrak r}_n)
 \end{equation} each consists of $q$-commuting elements. 
\end{Prop}

\proof This is proved in the same way as Proposition~{\ref{7.13}. \qed

\medskip

Our goal is to construct seeds out of these clusters using (and then 
generalizing) Proposition~\ref{seedth}.

With  Proposition~\ref{qcom} at hand, we are immediately given  the 
corresponding 
symplectic matrices 
\begin{equation}
{\mathcal L}_q^o({\mathfrak r}_1,\dots, {\mathfrak 
r}_{n-1},{\mathfrak r}_n)\textrm{ and }{\mathcal L}_q({\mathfrak r}_1,\dots, 
{\mathfrak 
r}_{n-1},{\mathfrak r}_n).
 \end{equation}
 The construction of the accompanying $B$-matrices
 \begin{equation}
{\mathcal B}_q^o({\mathfrak r}_1,\dots, {\mathfrak 
r}_{n-1},{\mathfrak r}_n)\textrm{ and }{\mathcal B}_q({\mathfrak r}_1,\dots, 
{\mathfrak 
r}_{n-1},{\mathfrak r}_n)
 \end{equation}
 takes a little more work, though in principle it is straightforward.
 The idea is in both cases to consider an element in the cluster as lying in a 
space
 \begin{eqnarray}
 {\mathcal C}_q^d({\mathfrak r}_i, {\mathfrak 
r}_{n-i},{\mathfrak r}_{n-i+1})\cup {\mathcal C}_q^u({\mathfrak r}_i, {\mathfrak 
r}_{i+1},{\mathfrak r}_{n-i})&\subseteq&{\mathcal C}_q({\mathfrak r}_i, 
{\mathfrak 
r}_{n-i},{\mathfrak r}_{n-i+1})\textrm{ or}\\
 {\mathcal C}_q^u({\mathfrak r}_i, {\mathfrak 
r}_{i+1},{\mathfrak r}_{n-i+1})\cup {\mathcal C}_q^d({\mathfrak r}_{i+1}, 
{\mathfrak 
r}_{n-i},{\mathfrak r}_{n-i+1})&\subseteq&{\mathcal C}_q^o({\mathfrak r}_i, 
{\mathfrak 
r}_{i+1},{\mathfrak r}_{n-i+1})\quad
 \end{eqnarray}
 as appropriate. Then we can use the corresponding matrices 
 ${\mathcal B}_q({\mathfrak r}_i, {\mathfrak 
r}_{n-i},{\mathfrak r}_{n-i+1})$ or ${\mathcal B}_q^o({\mathfrak r}_i, {\mathfrak 
r}_{i+1},{\mathfrak r}_{n-i+1})$ in the sense that one can extend these 
matrices 
to the full rank by inserting rows of zeros.
In this way, we can construct columns even for the troublesome  elements of the 
form  $E(a_{{\mathfrak r}_i}, a_{{\mathfrak r}_j})$ that may belong to such 
spaces. Indeed, we may start by including  $E(a_{{\mathfrak r}_{\frac{ n+0}2}}, 
a_{{\mathfrak r}_{\frac{n+2}2}})$ ($n$ even) or $E(a_{{\mathfrak 
r}_{\frac{n-1}2}}, a_{{\mathfrak r}_{\frac{n+1}2}})$ ($n$ odd)   in a such 
space 
in which they may be seen as mutable. Then these spaces have new non-mutable 
elements which can be handled by viewing them in appropriate spaces. The only 
ones which we cannot capture are the elements ${\det}_{s}^{{\mathfrak 
r}_1,{\mathfrak r}_n}=E(s_{{\mathfrak r}_1}, 
s_{{\mathfrak r}_n})$.

\begin{Def}In both cases, the elements  ${\det}_{s}^{{\mathfrak r}_1,{\mathfrak r}_n}$, $s\in 
Im(\pi_{{\mathfrak r}_1})$ 
are the non-mutable elements. We let ${\mathcal N}_q({\mathfrak r}_1, {\mathfrak r}_n)$ denote the set of these.
\end{Def}

\medskip
\begin{Prop}
 \begin{equation}
{\mathcal Q}_q({\mathfrak r}_1,\dots, {\mathfrak 
r}_{n-1},{\mathfrak r}_n)\textrm{ and } 
{\mathcal Q}_q^o({\mathfrak r}_1,\dots, {\mathfrak 
r}_{n-1},{\mathfrak r}_n)
 \end{equation}
 are quantum seeds.
 \end{Prop}

\bigskip

\section{Mutations}

Here is the fundamental setup: Let $\omega^{\mathfrak a}, \omega^{\mathfrak 
b},\omega^{\mathfrak c}\in W^p$ satisfy

\begin{equation}{\mathfrak a}<{\mathfrak c}\textrm{ and }  {\mathfrak a}\leq 
{\mathfrak b}\leq{\mathfrak c}.\end{equation}

\medskip

\begin{Def}\label{7.1bis}
A root $\gamma\in\triangle^+({\mathfrak c})$ is an {\bf increasing-mutation 
site}  
of 
$\omega^{\mathfrak b}\in W^p$ (in reference to  $({\mathfrak a},{\mathfrak 
b},{\mathfrak c})$)  if  there exists a reduced form  of 
$\omega^{\mathfrak c}$ as 
\begin{equation}
\omega^{\mathfrak c}=\hat\omega\sigma_\gamma\omega^{\mathfrak b}.
\end{equation}
Let $W^p\ni\omega^{{\mathfrak b}'}=\sigma_\gamma\omega^{\mathfrak b}$. It 
follows 
that 
\begin{equation}\label{94}
\omega^{{\mathfrak b}'}=\omega^{\mathfrak b}\sigma_{\alpha_s}
\end{equation}
for a unique $s\in Im(\pi_{{\mathfrak b}'})$. Such a site will henceforth  be 
called an ${\mathfrak m}^+$ site.

We will further say that $\gamma$ is a {\bf decreasing-mutation site}, or 
${\mathfrak m}^-$ site (in reference to  $({\mathfrak a},{\mathfrak 
b},{\mathfrak c})$)
of $\omega^{{\mathfrak b}}\in W^p$ in case there exists a rewriting of 
$\omega^{{\mathfrak b}}$ as $\omega^{{\mathfrak 
b}}=\sigma_\gamma\omega^{{\mathfrak b}''}$ with ${\mathfrak a}\leq 
\omega^{{\mathfrak b}''}\in W^p$. Here, \begin{equation}
\omega^{{\mathfrak b}}=\omega^{{\mathfrak b}''}\sigma_{\alpha_s}
\end{equation}
for a unique $s\in Im(\pi_{{\mathfrak b}})$. We view such sites as places where 
replacements are possible and will use the notation

\begin{equation}\label{m+}{\mathfrak m}^{+}_{{\mathfrak a},{\mathfrak c}}:({\mathfrak 
a},{\mathfrak 
b},{\mathfrak c})\to ({\mathfrak a},{\mathfrak b}',{\mathfrak c}),\end{equation} 
and 
\begin{equation}{\mathfrak m}^-_{{\mathfrak a},{\mathfrak c}}:({\mathfrak 
a},{\mathfrak 
b},{\mathfrak c})\to ({\mathfrak a},{\mathfrak b}'',{\mathfrak c}),\end{equation}
respectively, for the replacements while at the same time defining what we mean 
by replacements. 

Notice that ${\mathfrak a}={\mathfrak b}$ and ${\mathfrak b}'={\mathfrak c}$ 
are 
allowed in the first while  ${\mathfrak b}={\mathfrak c}$ and ${\mathfrak 
b}''={\mathfrak a}$ are allowed in the second.
 
Furthermore,
$${\mathfrak m}_{{\mathfrak a},{\mathfrak c}}:({\mathfrak a},{\mathfrak 
b},{\mathfrak c})\to ({\mathfrak a},{\mathfrak b}_1,{\mathfrak c})$$ denotes 
the 
composition of any finite number of such maps ${\mathfrak m}^{\pm}_{{\mathfrak 
a},{\mathfrak c}}$ (in any order, subject to the limitations at any step 
stipulated above)

We will further extend the meaning of ${\mathfrak m}_{{\mathfrak a},{\mathfrak 
c}}$ also to include the replacements
$${\mathcal C}_q({\mathfrak a},{\mathfrak b},{\mathfrak c})\to {\mathcal 
C}_q({\mathfrak a},{\mathfrak b}_1,{\mathfrak c}),$$ and even 
$${\mathcal Q}_q({\mathfrak a},{\mathfrak b},{\mathfrak c})\to {\mathcal 
Q}_q({\mathfrak a},{\mathfrak b}_1,{\mathfrak c}).$$At the seed level, we will 
refer to the replacements as {\bf Schubert mutations}.

Similarly, we can define maps ${\mathfrak m}^{o,\pm}_{{\mathfrak a},{\mathfrak 
c}}$, and after that mutations as composites $${\mathfrak m}^{o}_{{\mathfrak 
a},{\mathfrak c}}:{\mathcal Q}_q^o({\mathfrak a},{\mathfrak b},{\mathfrak c})\to 
{\mathcal Q}_q^o({\mathfrak a},{\mathfrak b}_1,{\mathfrak c}).$$
\end{Def}

\bigskip

We need to define another kind of replacement: 
Consider\begin{equation}\label{maxim}
{\mathfrak a}<{\mathfrak b}_1<{\mathfrak b}<{\mathfrak c}.
\end{equation}

\begin{Def}We say that  $({\mathfrak a},{\mathfrak b},{\mathfrak 
c})$ is a {\bf d-splitting} of $({\mathfrak a},{\mathfrak c})$ if 
$${\mathcal C}_q({\mathfrak a},{\mathfrak b},{\mathfrak 
c})={\mathcal C}_q({\mathfrak a},{\mathfrak c}).$$ In this case we will also say 
that    $({\mathfrak a},{\mathfrak c})$ is a {\bf 
d-merger} of  $({\mathfrak a},{\mathfrak b},{\mathfrak c})$.
\end{Def}

\medskip

To make this more definitive, one might further assume that  ${\mathfrak b}$ is 
maximal amongst those satisfying  (\ref{maxim}), but we will not need to do  
this here.

\medskip

Similarly,

\begin{Def}We say that  $({\mathfrak a},{\mathfrak b},{\mathfrak 
c})$ is a {\bf u-splitting} of $({\mathfrak a},{\mathfrak c})$ if 
$${\mathcal C}_q^o({\mathfrak a},{\mathfrak b},{\mathfrak 
c})={\mathcal C}_q^o({\mathfrak a},{\mathfrak c}).$$ Similarly, we will in this 
case also say that    $({\mathfrak a},{\mathfrak c})$ is a 
{\bf u-merger} of  $({\mathfrak a},{\mathfrak b},{\mathfrak 
c})$.
\end{Def}

\medskip

Our next definition combines the two preceding:

\begin{Def}\label{def84}A Schubert creation replacement
$$a^+_{{\mathfrak a},{\mathfrak c}}:({\mathfrak a},{\mathfrak 
c})\rightarrow ({\mathfrak a},{\mathfrak b}_1,{\mathfrak c})$$
consists in a d-splitting 
$$({\mathfrak a},{\mathfrak c})\rightarrow ({\mathfrak 
a},{\mathfrak b},{\mathfrak c})$$ followed by a replacement
$m_{{\mathfrak a},{\mathfrak c}}$ applied to $({\mathfrak a},{\mathfrak 
b},{\mathfrak c})$.
A Schubert annihilation replacement
$$a^-_{{\mathfrak a},{\mathfrak c}}:({\mathfrak a},{\mathfrak 
b}_1,{\mathfrak c})\rightarrow ({\mathfrak a},{\mathfrak c})$$ is 
defined as the reverse process.

Schubert creation/annihilation mutations $a^{o,\pm}_{{\mathfrak 
a},{\mathfrak c}}$ are defined analogously;
$$a^{o,+}_{{\mathfrak a},{\mathfrak c}}:
{\mathcal Q}_q^o({\mathfrak a},{\mathfrak c})\to {\mathcal 
Q}_q^o({\mathfrak a},{\mathfrak b}_1,{\mathfrak c}),
$$and 
$$a^{o,-}_{{\mathfrak a},{\mathfrak c}}:
{\mathcal Q}_q^o({\mathfrak a},{\mathfrak b}_1,{\mathfrak c})\to 
{\mathcal Q}_q^o({\mathfrak a},{\mathfrak c}).
$$
We finally extend these Schubert creation/annihilation mutations   into (we 
could do it more generally, but do not need to do so here)
$${\mathcal Q}_q({\mathfrak r}_1,\dots, {\mathfrak 
r}_{n-1},{\mathfrak r}_n)\rightarrow {\mathcal Q}_q({\mathfrak r}_1,\dots, 
{\mathfrak 
r}_{n-2},\dots,{\mathfrak r}_{n\pm 1})$$
by inserting/removing an ${\mathfrak r}_x$ between
${\mathfrak r}_{\frac{n}2}$ and  ${\mathfrak r}_{\frac{n}2+1}$ ($n$ even) or 
between ${\mathfrak r}_{\frac{n-1}2}$ and  ${\mathfrak r}_{\frac{n+1}2}$ 
($n$ odd). Similar maps are defined for the spaces 
 ${\mathcal Q}_q^o({\mathfrak r}_1,\dots, {\mathfrak 
r}_{n-1},{\mathfrak r}_n)$.
\end{Def}

\bigskip

In the sequel, we will encounter expressions of the form $\check 
B(u,v,s)$;

\medskip
\begin{equation}\label{99}
\check B(u,v,s)=E_{u\sigma_s\Lambda_s,v\Lambda_s}^{-1}E_{
u\Lambda_s,v\sigma_s\Lambda_s}^{-1}\prod_{a_{ks<0}} 
E_{u\Lambda_k,v\Lambda_k}^{a_{ks}}\end{equation} where \begin{equation}
E(u\Lambda_s,v\Lambda_s)\check 
B(u,v,s)=q^{-2(\Lambda_s,\alpha_s)}\check 
B(u,v,s)E(u\Lambda_s,v\Lambda_s),
\end{equation}
and where $\check B(u,v,s)$ commutes with all other 
elements in a given cluster.

\begin{Def}
We say that $\check B(u,v,s)$ implies the change 
$$E_{u\Lambda_s,v\Lambda_s}\to E_{u\sigma_s\Lambda_s,v\sigma_s\Lambda_s}.$$
\end{Def}

We will only encounter such changes where the set with 
$E_{u\Lambda_s,v\Lambda_s}$ removed from the initial cluster, and 
$E_{u\sigma_s\Lambda_s,v\sigma_s\Lambda_s}$ added, again is a cluster. 

We further observe that a (column) vector with $-1$ at positions corresponding 
to  $E_{u\sigma_s\Lambda_s,v\Lambda_s}$ and $E_{u\Lambda_s,v\sigma_s\Lambda_s}$ 
and $
{a_{ks}}$ at each position corresponding to a $E_{u\Lambda_k,v\Lambda_k}$ with 
$a_{ks}<0$ has the property that the symplectic form of the original cluster, 
when applied to it, returns a 
vector whose only non-zero entry is $-2(\Lambda_s,\alpha_s)$ at the position 
corresponding to  $E_{u\Lambda_s,v\Lambda_s}$. Hence, this can be a column 
vector of the $B$ of a potential compatible pair.

Even more can be ascertained: It can be seen that the last two lines of 
Theorem~\ref{toric} precisely states that with a $B$ matrix like that, the following holds:

\begin{Prop}The change $E_{u\Lambda_s,v\Lambda_s}\to E_{u\sigma_s\Lambda_s,v\sigma_s\Lambda_s}$ implied by $\check B(u,v,s)$
is the result of a BFZ mutation.
\end{Prop}

\medskip 

\begin{Thm} The Schubert mutation $${\mathcal Q}_q({\mathfrak a},{\mathfrak 
b},{\mathfrak c})\rightarrow{\mathcal Q}_q({\mathfrak a},{\mathfrak 
b}',{\mathfrak c})$$ implied by a replacement  ${\mathfrak m}^{+}_{{\mathfrak a},{\mathfrak c}}$ as in (\ref{m+}) is the result of  series of BFZ mutations.
\end{Thm}

\proof The number $s$ is given by (\ref{94}) and remains fixed throughout. We do the replacement in a number of steps. We set ${\mathcal Q}_q({\mathfrak a},{\mathfrak 
b},{\mathfrak c})={\mathcal Q}_q({\mathfrak a},{\mathfrak 
b},{\mathfrak c})(0)$ and perform changes
\begin{eqnarray}
 &{\mathcal Q}_q({\mathfrak a},{\mathfrak 
b},{\mathfrak c})={\mathcal Q}_q({\mathfrak a},{\mathfrak 
b},{\mathfrak c})(0)\rightarrow \\\nonumber& {\mathcal Q}_q({\mathfrak a},{\mathfrak 
b},{\mathfrak c})(1)\rightarrow \dots\rightarrow {\mathcal Q}_q({\mathfrak a},{\mathfrak 
b},{\mathfrak c})(t_0)={\mathcal Q}_q({\mathfrak a},{\mathfrak 
b}',{\mathfrak c}).
\end{eqnarray}

We will below see that $t_0=s_{\mathfrak b}-s_{\mathfrak a}-1$. We set

\begin{equation}
 \textrm{If } 0\leq t\leq t_o: {\mathcal Q}_q({\mathfrak a},{\mathfrak 
b},{\mathfrak c})(t)=({\mathcal C}_q({\mathfrak a},{\mathfrak 
b},{\mathfrak c})(t),{\mathcal L}_q({\mathfrak a},{\mathfrak 
b},{\mathfrak c})(t),{\mathcal B}_q({\mathfrak a},{\mathfrak 
b},{\mathfrak c})(t)).
\end{equation}
The intermediate seeds ${\mathcal Q}_q({\mathfrak a},{\mathfrak 
b},{\mathfrak c})(t)$ with $0<t<t_0$ are not defined by strings $\tilde{\mathfrak a}\leq \tilde{\mathfrak b}\leq \tilde{\mathfrak c}$. At each $t$-level, only one column is replaced when passing from ${\mathcal B}_q({\mathfrak a},{\mathfrak 
b},{\mathfrak c})(t)$ to ${\mathcal B}_q({\mathfrak a},{\mathfrak 
b},{\mathfrak c})(t+1)$, and here (\ref{77}) is applied. Of course, the whole ${\mathcal B}$ matrix is given by (\ref{72}) and (\ref{75}) for a suitable seed. 

\smallskip

Specifically, using (\ref{77}) we introduce a family of expressions $\check B$ as in (\ref{99}) 

\medskip

\begin{eqnarray}\label{b-conv-bar-t}{B}^{{\mathfrak a},{\mathfrak 
b}(t), {\mathfrak c}}_{m^+}(s,t)=
E_{\omega(s_{\mathfrak a}+t+1)\Lambda_s,\omega^{{\mathfrak b}}\Lambda_s}^{-1} 
E_{\omega(s_{\mathfrak a}+t)\Lambda_s,\omega^{{\mathfrak b}'}\Lambda_s}^{-1} 
\prod  
E_{\omega(s,s_{\mathfrak a}+t+1)\Lambda_j,\omega^{{\mathfrak b}'}\Lambda_j} 
^{-a_{js}}\\\nonumber
=(E^u_{\mathfrak b}(s,s_{\mathfrak a}+t+1)E^u_{{\mathfrak b}'}(s,s_{\mathfrak 
a}+t))^{-1}\prod E^u_{{\mathfrak b}}(j,\overline p(j,s,s_{\mathfrak 
a}+t+1))^{-a_{js}},
\end{eqnarray}
implying the changes 
\begin{equation}
E^u_{{\mathfrak b}}(s,s_{\mathfrak a}+t)\rightarrow E^u_{{\mathfrak 
b}'}(s,s_{\mathfrak a}+t+1).
\end{equation}

If $\omega(s,s_{\mathfrak a}+t+1)=u_t\sigma_s$ and $v=\omega^{\mathfrak b}$ 
then 
this corresponds to
\begin{equation}
\left((u_t\sigma_s\Lambda_s,v\Lambda_s)(u_t\Lambda_s,
v\sigma_s\Lambda_s)(u\Lambda_j,v\Lambda_j)^{a_{js}}\right)^{-1}
\end{equation}

Here are then in details how the changes are performed:

\begin{eqnarray*}{\mathrm Step}(0):&&\\{\mathcal C}_q({\mathfrak a},{\mathfrak 
b},{\mathfrak 
c})\ni E^d_{{\mathfrak a}}(s,s_{\mathfrak b}+1)&\rightarrow& E^u_{{\mathfrak 
b}'}(s,s_{\mathfrak a})\in {\mathcal C}_q({\mathfrak a},{\mathfrak 
b}(0),{\mathfrak c}) \ (renaming),\\
{B}_q^{{\mathfrak a},{\mathfrak 
b}, {\mathfrak c}}(s,s_{\mathfrak a})&\rightarrow&{B}^{{\mathfrak a},{\mathfrak 
b}(0), {\mathfrak c}}_{m^+}(s,0)  \ (renaming),
\\{\mathcal L}_q({\mathfrak a},{\mathfrak b},{\mathfrak 
c})&\rightarrow& {\mathcal L}_q({\mathfrak a},{\mathfrak b}(0),{\mathfrak 
c}) \ (renaming),
\\{\mathrm Step}(1): && (implied\ by \ {B}^{{\mathfrak a},{\mathfrak 
b}(0), {\mathfrak c}}_{m^+}(s,0) ),
  \\{\mathcal C}_q({\mathfrak a},{\mathfrak b}(0),{\mathfrak c}) \ni 
E^u_{{\mathfrak b}}(s,s_{\mathfrak a})&\rightarrow& E^u_{{\mathfrak 
b}'}(s,s_{\mathfrak a}+1) \in {\mathcal C}_q^d({\mathfrak a},{\mathfrak 
b}(1),{\mathfrak c}) , \\{B}_q^{{\mathfrak a},{\mathfrak 
b}, {\mathfrak c}}(s,s_{\mathfrak a}+1)&
\rightarrow&{B}^{{\mathfrak a},{\mathfrak 
b}(1), {\mathfrak c}}_{m^+}(s,1) (by\ (\ref{77})),
\\{\mathcal L}_q({\mathfrak a},{\mathfrak b}(0),{\mathfrak 
c})&\rightarrow& {\mathcal L}_q({\mathfrak a},{\mathfrak b}(1),{\mathfrak 
c}) \ (implied),\\{\mathrm Step}(2): && (implied\ by \ {B}^{{\mathfrak 
a},{\mathfrak 
b}(1), {\mathfrak c}}_{m^+}(s,1) ),
\\
{\mathcal C}_q^d({\mathfrak a},{\mathfrak b}(1),{\mathfrak c}) \ni 
E^u_{{\mathfrak 
b}}(s,s_{\mathfrak a}+1)&\rightarrow& E^u_{{\mathfrak b}'}(s,s_{\mathfrak 
a}+2)\in {\mathcal C}_q^d({\mathfrak a},{\mathfrak b}(2),{\mathfrak c}),\\\vdots\
\\{\mathrm Step}(t+1): && (implied\ by \ {B}^{{\mathfrak a},{\mathfrak 
b}(t), {\mathfrak c}}_{m^+}(s,t) ),\\
{\mathcal C}_q^d({\mathfrak a},{\mathfrak b}(t),{\mathfrak c}) \ni 
E^u_{{\mathfrak 
b}}(s,s_{\mathfrak a}+t)&\rightarrow& E^u_{{\mathfrak b}'}(s,s_{\mathfrak 
a}+t+1)\in {\mathcal C}_q^d({\mathfrak a},{\mathfrak b}(t+1),{\mathfrak c}) , 
\\{B}_q^{{\mathfrak a},{\mathfrak 
b}, {\mathfrak c}}(s,s_{\mathfrak a}+t)&
\rightarrow&{B}^{{\mathfrak a},{\mathfrak 
b}(t), {\mathfrak c}}_{m^+}(s,t)  (by\ (\ref{77})),\\{\mathcal L}_q({\mathfrak a},{\mathfrak 
b}(t),{\mathfrak 
c})&\rightarrow& {\mathcal L}_q({\mathfrak a},{\mathfrak b}(t+1),{\mathfrak 
c}) \ (implied).
\end{eqnarray*}

\medskip

The last step is $t=s_{\mathfrak b}-s_{\mathfrak a}-1$. ${\mathfrak 
b}(0)={\mathfrak b}$, ${\mathfrak b}(s_{\mathfrak b}-s_{\mathfrak 
a}-1)={\mathfrak b}'$.

\bigskip

It is easy to see that all intermediate sets indeed are seeds.

\medskip

\medskip

What is missing now is to connect, via a change of basis transformation of the
compatible pair, with the ``$E,F$'' matrices of \cite{bz}.  Here we notice that both terms
\begin{equation}
(E^u_{\mathfrak b}(s,s_{\mathfrak a}+t+1)E^u_{{\mathfrak b}'}(s,s_{\mathfrak 
a}+t))^{-1}(E^u_{\mathfrak b}(s,s_{\mathfrak a}+t))^{-1}
\end{equation}
and 

\begin{equation}
\prod E^u_{{\mathfrak b}}(j,\overline p(j,s,s_{\mathfrak 
a}+t+1))^{-a_{js}}(E^u_{\mathfrak b}(s,s_{\mathfrak a}+t))^{-1}
\end{equation}

have the same $q$-commutators as $E^u_{{\mathfrak b}'}(s,s_{\mathfrak a}+t+1)$. 
The two possibilities correspond to the two signs in formulas (3.2) and (3.3) in \cite{bz}.

\medskip

Indeed, the linear transformation 
\begin{equation}E(t):E^u_{\mathfrak b}(s,s_{\mathfrak a}+t)\rightarrow
-E^u_{\mathfrak b}(s,s_{\mathfrak a}+t+1)-E^u_{{\mathfrak b}'}(s,s_{\mathfrak 
a}+t)-E^u_{\mathfrak b}(s,s_{\mathfrak a}+t)
\end{equation}
results in a change-of-basis on the level of forms:
\begin{eqnarray}{\mathcal L}_q({\mathfrak a},{\mathfrak 
b}(t),{\mathfrak 
c})\rightarrow{\mathcal L}_q({\mathfrak a},{\mathfrak b}(t+1),{\mathfrak 
c})&=&E^T(t){\mathcal L}_q({\mathfrak a},{\mathfrak 
b}(t),{\mathfrak 
c})E(t),\\\nonumber {\mathcal B}_{m^+}^{{\mathfrak a},{\mathfrak 
b}(t),{\mathfrak 
c}}(s,t)\rightarrow{\mathcal B}_{m^+}^{{\mathfrak a},{\mathfrak b}(t+1),{\mathfrak 
c}}(s,t+1)&=&E(t){\mathcal B}_{m^+}^{{\mathfrak a},{\mathfrak 
b}(t),{\mathfrak 
c}}(s,t)F(t),
\end{eqnarray}
where $F(t)$ is a truncated part of $E(t)^T$ (the restriction to the mutable 
elements).

With this, the proof is complete. \qed

\medskip

\begin{Thm}
Any ${\mathcal Q}_q({\mathfrak r}_1,\dots, {\mathfrak 
r}_{n-1},{\mathfrak r}_n)$ can be obtained from ${\mathcal Q}_q({\mathfrak 
e},{\mathfrak p})$ as a sub-seed and any ${\mathcal Q}_q^o({\mathfrak r}_1,\dots, 
{\mathfrak 
r}_{n-1},{\mathfrak r}_n)$ can be obtained from ${\mathcal Q}_q^o({\mathfrak 
e},{\mathfrak p})$ as a sub-seed through a series of Schubert creation and 
annihilation 
mutations. 
These mutations are, apart from the trivial actions of renaming, splitting, 
merging, or simple restrictions, composites of BFZ-mutations.
\end{Thm}

\proof Apart from mergers and splittings (Definition~\ref{def84}), the mutations are composites of mutations of the form ${\mathcal Q}_q({\mathfrak a},{\mathfrak 
b},{\mathfrak c})\to {\mathcal Q}_q({\mathfrak a},{\mathfrak 
b}',{\mathfrak c})$. \qed

\medskip

\begin{Cor}\label{cor8.7}The algebras ${\mathcal A}_q^{d,{\mathfrak a},{\mathfrak 
c}}$ and ${\mathcal A}_q^{u,{\mathfrak a},{\mathfrak c}}$ are mutation equivalent 
and indeed are equal. We denote henceforth this algebra by ${\mathcal 
A}^{{\mathfrak a},{\mathfrak 
c}}$. This is the quadratic algebra  generated by the elements $\beta_{c,d}$ 
with $c_{\mathfrak a}<d\leq c_{\mathfrak c}$. 
\end{Cor}

\medskip

We similarly denote the corresponding skew-field of fractions by ${\mathcal 
F}_q^{{\mathfrak a},{\mathfrak 
c}}$.

\bigskip

\section{Prime}

\begin{Def}
\begin{equation}{\det}_{s}^{{\mathfrak a},{\mathfrak c}}:=E_{\omega^{\mathfrak 
a}\Lambda_s,\omega^{\mathfrak c}\Lambda_s}.\end{equation}
\end{Def}

\begin{Thm}\label{8.6} The 2 sided ideal $I({\det}_{s}^{{\mathfrak 
a},{\mathfrak 
c}})$ in ${\mathcal A}_q({\mathfrak a},{\mathfrak c})$ generated
by the covariant and non-mutable element ${\det}_{s}^{{\mathfrak a},{\mathfrak 
c}}$ is \underline{prime} for each $s$. 
\end{Thm}

\proof Induction. The induction start is trivially satisfied. Let us then 
divide the induction step into two cases. First, let  $Z_\gamma$ be an 
annihilation-mutation site of 
$\omega^{\mathfrak c}$ such that $\omega^{\mathfrak 
c}=\sigma_\gamma\omega^{{\mathfrak c}_1}=\omega^{{\mathfrak 
c}_1}\sigma_{\alpha_s}$
with $\omega^{{\mathfrak c}_1}\in W^p$. We have clearly ${\mathcal 
A}_q({\mathfrak a},{\mathfrak c})= {\mathcal A}_q({\mathfrak a},{\mathfrak 
c}_1)\cup I({\det}_{s}^{{\mathfrak a},{\mathfrak c}})$. Furthermore, ${\mathcal 
A}_q({\mathfrak a},{\mathfrak c})\setminus {\mathcal 
A}_q({\mathfrak a},{\mathfrak c}_1) =I_\ell(Z_\gamma)$, where 
$I_\ell(Z_\gamma)$ denotes the left ideal generated by $Z_\gamma$. We might as 
well consider the right ideal,
but not the 2-sided ideal since in general there will be terms ${\mathcal R}$ of lower order, c.f. Theorem~\ref{4.1}.

It follows that \begin{equation}\label{Z}{\det}_{s}^{{\mathfrak a},{\mathfrak 
c}}=M_1Z_\gamma 
+M_2\end{equation} where $M_1,M_2\in {\mathcal 
A}_q({\mathfrak a},{\mathfrak c}_1)$ and 
$M_1\neq0$.
Indeed, $M_1$ is a non-zero multiple of ${\det}_{s}^{{\mathfrak a},{\mathfrak 
c}_1}$. (If
$s_{\mathfrak c}=1$ then $M_1=1$ and $M_2=0$.) We also record, partly for 
later use, 
that $Z_\gamma$ $q$-commutes with everything up to correction terms from 
${\mathcal 
A}_q({\mathfrak a},{\mathfrak c}_1)$.

\medskip

Notice that we use Corollary~\ref{cor8.7}.
\medskip

Now consider an equation
\begin{equation}
 {\det}_{s}^{{\mathfrak a},{\mathfrak c}}p_1=p_2p_3
\end{equation}
with $p_1,p_2,p_3\in {\mathcal A}_q({\mathfrak a},{\mathfrak 
c})$. Use (\ref{Z}) to write for each $i=1,2,3$
\begin{equation}
 p_i=\sum_{k=0}^{n_i}({\det}_{s}^{{\mathfrak a},{\mathfrak c}})^kN_{i,k} 
\end{equation}
where each $N_{i,k}\in {\mathbf L}_q({\mathfrak a},{\mathfrak 
c}_1)$ and assume that $N_{i,0}\neq0\textrm{ for }i=2,3$
Then $0\neq N_{0,2}N_{0,3}\in {\mathbf L}_q({\mathfrak a},{\mathfrak c}_1)$. At the same time, 
\begin{equation}
 N_{0,2}N_{0,3}=\sum_{k=1}^{n_i}({\det}_{s}^{{\mathfrak a},{\mathfrak 
c}})^k\tilde N_{i,k} 
\end{equation}
for certain elements $\tilde N_{i,k} \in {\mathbf L}_q({\mathfrak a},{\mathfrak 
c}_1)$.

Using the linear independence (\cite[Proposition~10.8]{bz}) we easily get a 
contradiction by looking at the leading term in ${\det}_{s}^{{\mathfrak 
a},{\mathfrak 
c}})$.

\smallskip

Now in the general case, the $s$ in ${\det}_{s}^{{\mathfrak a},{\mathfrak c}}$  is given 
and we may write ${\omega}^{\mathfrak c}={\omega}^{{\mathfrak 
c}_2}\sigma_{s}\tilde{\omega}$ where $\sigma_s$ does not occur in 
$\tilde{\omega}$. Let ${\omega}^{{\mathfrak c}_1}={\omega}^{{\mathfrak 
c}_2}\sigma_s$. It is clear that ${\det}_{s}^{{\mathfrak a},{\mathfrak 
c}}={\det}_{s}^{{\mathfrak a},{\mathfrak c}_1}$ and by the previous, 
${\det}_{s}^{{\mathfrak a},{\mathfrak c}_1}$ is prime in 
${\mathcal A}_q({\mathfrak a},{\mathfrak c}_1)$. We have that ${\mathcal 
A}_q({\mathfrak a},{\mathfrak c}_1)$ is 
an algebra in its own right. Furthermore,
\begin{equation}{\mathcal A}_q({\mathfrak a},{\mathfrak c})={\mathcal 
A}_q({\mathfrak a},{\mathfrak c}_1)[Z_{\gamma_1},\dots,Z_{\gamma_n}],
\end{equation}
where the Lusztig elements $Z_{\gamma_1},\dots,Z_{\gamma_n}$  are bigger than 
the generators of ${\mathcal A}_q({\mathfrak a},{\mathfrak c}_1)$. In a  PBW basis 
we can put them to the right. They even generate a quadratic algebra 
$\tilde{\mathcal A}_q$ in their own right! The equation we need to consider are 
of 
the form
\begin{equation}p_1p_2={\det}_{s}^{{\mathfrak a},{\mathfrak c}_1}p_3
\end{equation} 
with $p_1,p_2,p_3\in {\mathcal A}_q({\mathfrak a},{\mathfrak c})$. The claim that 
at least one 
of $p_1,p_2$ contains a factor of ${\det}^{{\mathfrak r}_1}_{q,s}$ follows by 
easy induction on the $\tilde{\mathcal A}_q$ degree of $p_1p_2$, i.e. the sum of 
the $\tilde{\mathcal A}_q$ degrees of $p_1$ and $p_2$.     \qed

\bigskip

\section{Upper}

Let $\omega^{\mathfrak a}, \omega^{\mathfrak c}\in W^p$ and ${\mathfrak 
a}<{\mathfrak 
c}$.

\begin{Def}
The  cluster algebra ${\mathbf A}_q({\mathfrak a},{\mathfrak c})$ is the 
${\mathbb 
Z}[q]$-algebra generated in the 
space ${\mathcal F}_q({\mathfrak a},{\mathfrak c})$ by the inverses of the 
non-mutable elements 
${\mathcal N}_q({\mathfrak a},{\mathfrak c})$ together with the union of the sets 
of  all  
variables
obtainable from the initial seed ${\mathcal Q}_q({\mathfrak a},{\mathfrak c})$ by 
composites of 
quantum
Schubert mutations. 
(Appropriately applied)
\end{Def}

\medskip
Observe that we include ${\mathcal N}_q({\mathfrak a},{\mathfrak c})$ in the set 
of variables.
\medskip

\begin{Def}
The upper cluster algebra ${\mathbf U}_q({\mathfrak a},{\mathfrak c})$ connected 
with 
the same pair $\omega^{\mathfrak a}, \omega^{\mathfrak c}\in W^p$ is the 
${\mathbb Z}[q]$-algebra in  ${\mathcal F}_q({\mathfrak a},{\mathfrak c})$ given 
as 
the intersection of  all the Laurent algebras of the 
sets of variables
obtainable from the initial seed ${\mathcal Q}_q({\mathfrak a},{\mathfrak c})$ by 
composites of 
quantum Schubert mutations. (Appropriately applied)
\end{Def}

\medskip

\begin{Prop}
$${\mathcal A}_q({\mathfrak a},{\mathfrak c})\subseteq {\mathbf A}_q({\mathfrak 
a},{\mathfrak c})\subset {\mathbf U}_q({\mathfrak a},{\mathfrak c}).$$
\end{Prop}

\proof The first inclusion follows from \cite{leclerc}, the other is the quantum 
Laurent phenomenon. \qed

\medskip

\begin{Rem}
Our terminology may seem a bit unfortunate since the notions of a cluster 
algebra and an
upper cluster algebra already have been introduced by Berenstein and Zelevinsky 
in
terms of all mutations. We only use quantum line mutations which form a proper
subset of the set of all quantum mutations. However, it will be a corollary to 
what
follows that the two notions in fact coincide, and for this reason we
do not introduce some auxiliary notation. 
\end{Rem}

\medskip

\begin{Thm}
 $${\mathbf U}_q({\mathfrak a},{\mathfrak c})={\mathcal A}_q({\mathfrak 
a},{\mathfrak c})[({\det}_{s}^{{\mathfrak a},{\mathfrak c}})^{-1}; s\in 
Im(\pi_{{\mathfrak c}})].$$
 \end{Thm}

\proof This follows by induction on $\ell(\omega^{\mathfrak c})$ (with start at 
$\ell(\omega^{\mathfrak a})+1$)  in the same way 
as in the proof of \cite[Theorem~8.5]{jz}, but for clarity we give the details: 
Let the notation and assumptions be as in the proof of Theorem~\ref{8.6}. First 
of all, the induction start is trivial since we there are looking at the 
generator of a Laurent quasi-polynomial algebra. Let then $u\in {\mathbf 
U}_q({\mathfrak a},{\mathfrak c})$. We will argue by contradiction, and just as 
in 
the proof 
of \cite[Theorem~8.5]{jz}, one readily sees that one may assume that $u\in 
{\mathcal  A}_q({\mathfrak a},{\mathfrak c}_1)[({\det}_{s}^{{\mathfrak 
a},{\mathfrak c}_1})^{-1}, {\det}_{s}^{{\mathfrak a},{\mathfrak c}}]$. Using 
(\ref{Z}) we may now write
\begin{equation} \label{neg1}
u=\left(\sum_{i=0}^K Z_\gamma^ip_i({\det}_{s}^{{\mathfrak a},{\mathfrak 
c}_1})^{k_i}\right)({\det}_{s}^{{\mathfrak a},{\mathfrak c}_1})^{-\rho},
\end{equation}
with $p_i\in {\mathcal A}_q({\mathfrak a},{\mathfrak c}_1)$, $p_i\notin 
I({\det}_{s}^{{\mathfrak a},{\mathfrak c}_1})$, and $k_i\geq0$. Our assumption 
is that $\rho>0$. recall that the elements
${\det}_{s}^{{\mathfrak a},{\mathfrak c}_1}$ and ${\det}_{s}^{{\mathfrak 
a},{\mathfrak c}}$ are
covariant and define prime ideals in the appropriate algebras.
Using the fact that ${\mathbf U}_q({\mathfrak a},{\mathfrak c})$ is an algebra 
containing ${\mathcal A}_q({\mathfrak a},{\mathfrak c})$, we can assume that the 
expression in the left bracket in (\ref{neg1}) is not in
$I({\det}_{s}^{{\mathfrak a},{\mathfrak c}})$ and we may further assume that 
$p_i\neq0\Rightarrow
k_i<\rho$. To wit, one can remove the factors of ${\det}_{s}^{{\mathfrak 
a},{\mathfrak c}}$, then
remove the terms with $k_i\geq \rho$, then possibly repeat this process a 
number 
of
times.

Consider now the cluster ${\mathcal C}_q^{u}({\mathfrak a},{\mathfrak c})$. We 
know that $u$ can be written as a Laurent quasi-polynomial in the elements 
of 
${\mathcal C}_q^{u}({\mathfrak a},{\mathfrak c})$. By factoring out, we can 
then write

\begin{equation}\label{neg2}
u=p\prod_{(c,d)\in{\mathbb 
U}^{u,{\mathfrak a},{\mathfrak c}}}(E^u_{\mathfrak c}(c,d))^{-\alpha_{c,d}},
\end{equation}
where $p\in {\mathcal A}_q({\mathfrak a},{\mathfrak c})$, and 
$\alpha_{c,d}\geq0$. 
We will 
compare this to (\ref{neg1}). For the sake of this argument set $\tilde{\mathbb 
U}^{u,{\mathfrak a},{\mathfrak c}})=\{(c,d)\in {\mathbb 
U}^{u,{\mathfrak a},{\mathfrak c}})\mid \alpha_{c,d}>0\}$. 

Of course,  ${\det}_s^{{\mathfrak a},{\mathfrak c}}\in {\mathcal 
C}^{u}({\mathfrak e},{\mathfrak r})$.

``Multiplying across'', we get from (\ref{neg1}) and (\ref{neg2}), absorbing 
possibly some terms into $p$:
\begin{equation}
(\sum_{i=0}^K Z^ip_i({\det}^{{\mathfrak a},{\mathfrak 
c}_1}_{s})^{k_i})\prod_{(c,d)\in\tilde{\mathbb 
U}^{u,{\mathfrak a},{\mathfrak c}}}(E^u_{\mathfrak 
c}(c,d))^{\alpha_{c,d}}=p({\det}^{{\mathfrak a},{\mathfrak c}_1}_{s})^{\rho}.
\end{equation}
Any factor of ${\det}^{{\mathfrak a},{\mathfrak c}}_{s}$ in $p$ will have to be 
canceled by a similar
factor of $E^u_{\mathfrak c}(s,0)$ in the left-hand side, so we can assume that 
$p$ does not contain no
factor of ${\det}^{{\mathfrak a},{\mathfrak c}}_{s}$. After that we can assume 
that $(s,0)\notin
\tilde{\mathbb 
U}^{u,{\mathfrak a},{\mathfrak c}}$ since clearly ${\det}^{{\mathfrak 
a},{\mathfrak c}_1}_{s}\notin
I({\det}^{{\mathfrak a},{\mathfrak c}}_{s})$. Using that $k_i<\rho$ it follows 
that there must be a
factor of $({\det}^{{\mathfrak a},{\mathfrak c}_1}_{s})$ in
$\prod_{(c,d)\in\tilde{\mathbb 
U}^{u,{\mathfrak a},{\mathfrak c}}}(E^u_{\mathfrak c}(c,d))^{\alpha_{c,d}}$. 
Here, as but noticed, $d=0$ is excluded. The other 
terms do not contain $Z_{s,1}$ but $({\det}^{{\mathfrak a},{\mathfrak 
c}_1}_{s})$ does. This is an obvious contradiction. \qed

\bigskip

\section{The diagonal of a quantized minor}

\begin{Def}Let ${\mathfrak a}<{\mathfrak b}$. The diagonal, ${\mathbb 
D}_{\omega^{\mathfrak 
a}(\Lambda_s),\omega^{\mathfrak b}(\Lambda_s)}$, of 
$E_{\omega^{\mathfrak a}(\Lambda_s),\omega^{\mathfrak 
b}(\Lambda_s)}$ is set to 
\begin{equation}
{\mathbb D}_{\omega^{\mathfrak a}(\Lambda_s),\omega^{\mathfrak 
b}(\Lambda_s)}=q^{\alpha}Z_{s,s_{\mathfrak a}+1}\cdots Z_{s,s_{\mathfrak b}},
\end{equation}
where 
\begin{equation}
Z_{s,s_{\mathfrak b}}\cdots Z_{s,s_{\mathfrak 
a}+1}=q^{2\alpha}Z_{s,s_{\mathfrak a}+1}\cdots Z_{s,s_{\mathfrak b}} + {\mathcal R}
\end{equation}
where the terms ${\mathcal R}$ are of lower order
\end{Def}

\medskip

\begin{Prop}
$$E_{\omega^{\mathfrak a}(\Lambda_s),\omega^{\mathfrak 
b}(\Lambda_s)}={\mathbb D}_{\omega^{\mathfrak a}(\Lambda_s),\omega^{\mathfrak 
b}(\Lambda_s)}+{\mathcal R}$$
The terms in ${\mathcal R}$ are  of lower order in our ordering induced by 
$\leq_L$. They can in theory be determined from the 
fact
that the full
polynomial belongs to the dual canonical
basis. (\cite{bz},\cite{leclerc}).
\end{Prop}

\medskip

\proof   We prove this by induction on the length $s_{\mathfrak b}-s_{\mathfrak 
a}$ of any $s$-diagonal. When this length is $1$ we have at most a quasi-polynomial algebra and here the case is clear. Consider then a creation-mutation site where we go from length $r$ to $r+1$: Obviously,  it is only  the very last determinant we need to consider. Here we use the equation in 
Theorem~\ref{3.2} but reformulate it in terms of  the elements $E_{\xi,\eta}$, cf. 
Theorem~\ref{toric}. 

Set $\omega^{{\mathfrak b}_1}=\omega^{{\mathfrak b}}\sigma_s$ and consider 
$E_{\omega^{\mathfrak a}(\Lambda_s),\omega^{{\mathfrak 
b}_1}(\Lambda_s)}$. Its weight is
given as
$$\omega^{{\mathfrak 
b}_1}(\Lambda_s)-\omega^{\mathfrak a}(\Lambda_s)=\beta_{s,s_{\mathfrak a}+1} 
\dots+\beta_{s,s_{\mathfrak b}+1}.$$
In the recast version of Theorem~\ref{3.2}, the terms on the left hand 
side are covered by the induction hypothesis. The second term on the right hand 
side contains no element of the form $Z_{s,s_{{\mathfrak b}_1}}$ and it follows 
that we have an equation

\begin{equation}(Z_{s,s_{\mathfrak a}+2}\cdots Z_{s,s_{\mathfrak 
b}})E_{\omega^{\mathfrak a}(\Lambda_s),\omega^{{\mathfrak 
b}_1}(\Lambda_s)}=(Z_{s,s_{\mathfrak a}+2}\cdots Z_{s,s_{\mathfrak b}+1}) 
(Z_{s,s_{\mathfrak a}+1}\cdots Z_{s,s_{\mathfrak b}})
 +{\mathcal R}.
\end{equation}
The claim follows easily from that. \qed

\medskip

Recall that in the associated quasi-polynomial algebra is the algebra with 
relations corresponding to the top terms, i.e., colloquially  speaking, setting 
the lower order terms ${\mathcal R}$ equal to $0$. Let
\begin{equation}
 {d}_{\omega^{\mathfrak r}_{s,t_1}(\Lambda_s),\omega^{\mathfrak 
r}_{s,t}(\Lambda_s)}=z_{s,t_1+1}\cdots z_{s,t}.
\end{equation}

\medskip

The following shows the importance of the diagonals:

\smallskip

\begin{Thm} 
\begin{eqnarray}
d_{u\cdot\Lambda_{i_0},v\cdot\Lambda_{i_0}}d_{u_1\cdot\Lambda_{i_1},
v_1\cdot\Lambda_{i_1}}&=&
q^G d_{u_1\cdot\Lambda_{i_1},v_1\cdot\Lambda_{i_1}} 
d_{u\cdot\Lambda_{i_0},v\cdot\Lambda_{i_0}}\Leftrightarrow\\
{\mathbb D}_{u\cdot\Lambda_{i_0},v\cdot\Lambda_{i_0}}{\mathbb
D}_{u_1\cdot\Lambda_{i_1},v_1\cdot\Lambda_{i_1}}&=&
q^G {\mathbb D}_{u_1\cdot\Lambda_{i_1},v_1\cdot\Lambda_{i_1}}
{\mathbb
D}_{u\cdot\Lambda_{i_0},v\cdot\Lambda_{i_0}}+ {\mathcal R}
 \end{eqnarray}In particular, if the two elements  
${E}_{u\cdot\Lambda_{i_0},v\cdot\Lambda_{i_0}}{E}_{u_1\cdot\Lambda_{i_1},
v_1\cdot\Lambda_{i_1}}$ $q$-commute:
\begin{equation}
{E}_{u\cdot\Lambda_{i_0},v\cdot\Lambda_{i_0}}{E}_{u_1\cdot\Lambda_{i_1},
v_1\cdot\Lambda_{i_1}}=
q^G {E}_{u_1\cdot\Lambda_{i_1},v_1\cdot\Lambda_{i_1}}
{E}_{u\cdot\Lambda_{i_0},v\cdot\Lambda_{i_0}}
 \end{equation}then $G$ can be computed in the associated quasi-polynomial 
algebra:
 \begin{equation}
d_{u\cdot\Lambda_{i_0},v\cdot\Lambda_{i_0}}d_{u_1\cdot\Lambda_{i_1},
v_1\cdot\Lambda_{i_1}}=
q^G D_{u_1\cdot\Lambda_{i_1},v_1\cdot\Lambda_{i_1}} 
d_{u\cdot\Lambda_{i_0},v\cdot\Lambda_{i_0}}.\end{equation}
\end{Thm}

\begin{Rem}
One can also compute $G$ directly using the formulas in \cite{bz}.
\end{Rem}

\medskip

\begin{Rem}
The elements $E_{\xi,\eta}$ that we consider belong to the dual canonical 
basis. 
As such, they can in principle be determined from the highest order terms 
${\mathbb D}_{\xi,\eta}$.
\end{Rem}

\bigskip

\bigskip

\section{Litterature}

\end{document}